\input amstex
\documentstyle{amsppt}
\pagewidth{6.4in}
\vsize8.5in
\parindent=6mm
\parskip=3pt
\baselineskip=14pt
\tolerance=10000
\hbadness=500
\NoRunningHeads
\loadbold
\topmatter
\title
Endpoint mapping properties of spherical maximal operators
\endtitle
\author Andreas Seeger \ \ \ Terence Tao \ \ \ James Wright\endauthor
\thanks
The first author is supported 
in part by a grant from the National Science
Foundation.  The second author is a Clay Prize fellow
and is supported by the Sloan and Packard foundations. 
  \endthanks
\address
Department of Mathematics,
University of Wisconsin, Madison, WI 53706-1388
\endaddress
\email seeger\@math.wisc.edu\endemail
\address
Department of Mathematics,
University of California, Los Angeles, CA  90095-1555, USA
\endaddress
\email tao\@math.ucla.edu\endemail
\address
Department of Mathematics and Statistics,
University of Edinburgh, King's Building, Mayfield Road, Edinburgh EH3 9JZ, U.K.
\endaddress
\email jimw\@maths.unsw.edu.au\endemail
\endtopmatter
\document

\def\s{\sigma}

\define\cf{{\it cf}}

\define\dist{{\text{\rm dist}}}

\define\rad{{\text{\rm rad}}}

\define\inn#1#2{\langle#1,#2\rangle}

\define\lcontr{\rfloor}
\define\lco#1#2{{#1}\lcontr{#2}}
\define\lcoi#1#2{\imath({#1}){#2}}
\define\rco#1#2{{#1}\rcontr{#2}}

\define\bin#1#2{{\pmatrix {#1}\\{#2}\endpmatrix}}
\define\meas{{\text{\rm meas}}}

\define\card{\text{\rm card}}
\define\lc{\lesssim}
\define\gc{\gtrsim}


\define\eps{\varepsilon}

\define\la{\lambda}

\define\om{\omega}

\define\fI{{\frak I}}

\define\fM{{\frak M}}
\define\fN{{\frak N}}

\define\fR{{\frak R}}

\define\fW{{\frak W}}


\define\bbR{{\Bbb R}}

\define\bbZ{{\Bbb Z}}

\define\cA{{\Cal A}}
\define\cB{{\Cal B}}
\define\cC{{\Cal C}}
\define\cD{{\Cal D}}
\define\cE{{\Cal E}}

\define\cG{{\Cal G}}

\define\cJ{{\Cal J}}

\define\cM{{\Cal M}}
\define\cN{{\Cal N}}

\define\cQ{{\Cal Q}}
\define\cR{{\Cal R}}

\define\cV{{\Cal V}}
\define\cW{{\Cal W}}




\define\Ga{\Gamma}

\head{\bf 1. Introduction}\endhead
For a function $f\in L^p(\Bbb R^d)$, $d\ge 2$, we consider
the spherical means
$$
\cA_t f(x)\,=\,\int_{S^{d-1}}f(x-ty) d\sigma(y)
\tag 1.1
$$
where $d\sigma$ is the rotationally invariant measure on $S^{d-1}$, normalized
such that $\sigma(S^{d-1})=1$.
We wish to study the  question of
pointwise convergence as $t\to 0$
where the radii $t$ are restricted to a subset $E$ of $(0,\infty)$.
Pointwise convergence is established from boundedness properties of the maximal function
$$
\Cal M_E f(x)\,=\,\sup_{t\in E} |\cA_t f(x)|
$$
for $f\in L^p(\Bbb R^d)$.

Stein \cite{14} showed that for  $E=\Bbb R_+$
the maximal operator $\Cal M_E$ is
bounded  on $L^p$ if and only if
 $p>d/(d-1)$, $d\ge 3$; the same result for the case
$d=2$ was later proved by Bourgain \cite{2}.
The critical exponent $p(E)$ for $L^p$ boundedness of $\Cal M_E$,
  for any set $E\in(0,\infty)$, was determined  by 
Seeger, Wainger and Wright \cite{12}.
It is computed using a dilation invariant
notion of Minkowski-dimension. In order to describe the result  we
 let $N(E,\delta)$ be the
$\delta$-entropy number of $E$, that is the minimal number of intervals of
length $\delta$
needed to cover $E$ (we shall always redefine
 $N(\emptyset,\delta)=1$).
Define $$E^k=[2^{k}, 2^{k+1})\cap E$$
and
$$p(E)\,=\, 1+\frac 1{d-1}\Bigl(
\sup_{\delta> 0}
\sup_{k\in\Bbb Z} \frac
{\log N(E^k,2^k\delta)}{\log \delta^{-1}}
\Bigr).
$$
Then $\Cal M_E$ is bounded on $L^p$ for $p>p(E)$ and unbounded on $L^p$ if
$p<p(E)$. Moreover  various $L^p$ results were proven  in \cite{12}
for the critical exponent $p=p(E)$; however these results  fell short
of being
 necessary and sufficient.

For the case
that our maximal operator acts only on {\it radial}
functions sharp endpoint estimates in almost all cases  have been 
obtained in \cite{13}.
The relevant condition for $1<p<d/(d-1)$ turned out to be
\proclaim{Condition $(\cC_{p,q})$}
$$\alignat2
&\sup_j\Big(\sum_{n\ge 0}[N(E^{j+n},2^j)]^{q/p} 
2^{-n(d-1)q/p'}\Big)^{1/q}<\infty &&\quad\text{ if $p\le q< \infty,$}
\tag 1.2 \\
&\sup\Sb k\in\Bbb Z\\ \delta> 0\endSb
 N(E^k,2^{k}\delta)^{1/p}\delta^{(d-1)/p'}<\infty
&&\quad \text{ if $q= \infty.$}
\tag 1.3
\endalignat
$$
\endproclaim

It is shown  in \cite{13} that
 for $\cM_E$ to map $L^p_\rad$ to the Lorentz space $L^{p,q}$,
$1<p<d/(d-1)$, $p\le q\le \infty$
it is necessary and sufficient that condition $(\cC_{p,q})$ holds.
The necessity can be shown by testing $\cM_E$ on characteristic functions of small balls.
 Observe that
($\cC_{p,\infty}$) is the limiting case of
($\cC_{p,q}$)
 as $q\to\infty$.
For $p=d/(d-1)$ there are different characterizations for $L^p_\rad\to L^{p,q}$
boundedness,
at least when $d>2$.

The main purpose of this paper is to   
prove analogues of the  $L^p_\rad\to L^{p}$ and
$L^p_\rad\to L^{p,\infty}$ endpoint estimates
for general functions in $L^p$, assuming however
an additional regularity assumption (see hypothesis ($\cR_{p}$) below).
The main general results for $1<p\le d/(d-1)$  are stated in
Theorem I, II, III and IV below.
The case  where each set $E^k=E\cap[2^k, 2^{k+1}]$ is a convex sequence
serves as a model case (see \S8 below). In particular we have

\proclaim{Theorem 1.1 }

(i) Let $0<\alpha<\infty$  and let
$$E(\alpha)=\{2^k(1+\nu^{-\alpha}): k\in \Bbb Z,\nu\in \Bbb Z^+\}.
\tag 1.4$$ 
Then $M_{E(\alpha)}$ is of weak type $(p,p)$
if and only if  $p\ge 1+[(d-1)(\alpha+1)]^{-1}$.

(ii) Let $1/(d-1)< \beta<\infty$  and let
$$\widetilde E(\beta) =\{2^k (1+\log^{-\beta}(2+\nu)): 
k\in \Bbb Z,\nu\in \Bbb Z^+\} .
\tag 1.5 $$  Then  $M_{\widetilde E(\beta)}$ is of weak type
$(p,p)$ if and only if $p\ge d/(d-1)$.
\endproclaim

\remark{Remarks}

(a) 
Only the endpoint cases $p= 1+[(d-1)(\alpha+1)]^{-1}$ and $p=d/(d-1)$
are new. When $\beta < 1/(d-1)$, 
 $M_{\widetilde E(\beta)}$  fails to be weak type
$(d/(d-1), d/(d-1))$. The case $\beta = 1/(d-1)$ remains open.

(b)   For $p=1$ a slight variant  
was obtained by M. Christ who proved that the lacunary spherical maximal 
operator (with $E=\{2^k:k\in \Bbb Z\}$) maps the Hardy space 
$H^1$ to $L^{1,\infty}$. This can be deduced from a simple
 modification  of the proof below, and in fact the weak type estimates in \S5 are extensions of Christ's argument.

(c) It is not   known whether the lacunary spherical maximal function maps $L^1$ to $L^{1,\infty}$. The closest known 
result is a weak type $L\log\log L$ inequality proved  by the authors 
in \cite{11}.

\endremark

We shall now formulate a technical result on $L^p$ 
boundedness for $\cM_E$ which is only a minor improvement 
of the result in \cite{12}. It gives a reasonably sharp 
but not yet definitive 
estimate
for general sets $E$ of dilations. 
It will be applied however to sets which tend to be  
much thinner than the original sets.

\proclaim{Proposition 1.2}
Suppose that $d\ge 2$ and  $1<p\le d/(d-1)$.
Suppose that $\{\omega_j\}_{j=0}^\infty$ is a sequence of positive numbers
satisfying
$$\sum_{j\ge 0}\omega_j^{-{p'}}\le 1
\tag 1.6$$
and suppose that
$$
\sup_{k\in\bbZ} \sum_{j\ge 0} \omega_j^p
N(E^{k+j}, 2^{k}) 2^{-j(d-1)p/p'} \le A_0^p .
\tag 1.7
$$
Then $\cM_E$ is bounded on  $L^p(\Bbb R^d)$, with
 operator norm dominated by
$CA_0$.
\endproclaim

We now  describe our regularity assumption and begin with  the following

\definition{Definition} (i) A set $J\subset \Bbb R^+$  is  {\it equally spaced}
with  width $\delta$ and possible deviation $C > 1$   if
for  all $t\in J$  the inequalities
$$
C^{-1} \delta \le \dist(t, J\setminus\{t\})
\le
C \delta
\tag 1.8
$$
hold.

(ii) A family $\cJ=\{J\}$  of subsets of $\bbR^+$ is {\it uniformly equally spaced}
if  for every $J\in \cJ$ there is a $\delta=\delta(J)>0$  so that (1.8) holds with $\delta(J)$ 
and  a
 constant $C$
independent of $J$.

(iii) Let $J$ be an equally spaced subset of $\Bbb R^+$.
Then we call $a_J=\inf J$ and $b_J=\sup J$ the endpoints of $J$.

(iv) Let $\cJ$ be uniformly equally spaced family of subsets of $\Bbb R^+$.  Then we denote by
$\cD(\cJ)$ the set of endpoints $\cD(\cJ)=\cup_{J\in\cJ}\{a_J,b_J\}$.
\enddefinition

Our regularity assumption will say that
 each
$E^k$  can be split into ``not too many'' equally spaced  sets.
This gives a large class of examples, since in general the sets
$\cD^k$  of endpoints  are often much thinner than the sets $E^k$.

\proclaim{ Regularity hypothesis $(\cR_{p})$ }

 $E$ satisfies hypothesis $(\cR_p)$ if
for each $k$ there is a collection $\cJ^k=\{J\}$  of subsets of $[2^k,2^{k+1}]$ such that
$E^k\subset \cup_{J\in\cJ^k}J$ and
the following three  conditions are satisfied.

\roster

\item"{{\it (a)}}"  The family $\{J: J \in \cup_{k\in \bbZ}\cJ^k\}$ is uniformly equally spaced
(with uniform possible deviation $C$).

\item"{{\it (b)}}"
There is  a positive sequence $\om=\{\om_j\}_{j=0}^\infty$  with
$\sum_{j=0}^\infty\om_j^{-p'}\le 1$
so that
the sets of endpoints
$\cD^k\equiv \cD(\cJ^k)=\cup_{J\in\cJ^k}\{a_J,b_J\}$ satisfy
$$\sup_{k\ge 0} \Big(\sum_{j\ge 0} \big[ N(\cD^{k+j}, 2^{k})
2^{-j(d-1)p/p'} \omega_j^p\Big)^{1/p} \le C_0<\infty.
\tag 1.9
$$

\item"{{\it (c)}}"
Let $\cJ^k_\mu$ denote the subfamily of all $J\in \cJ^k$ which are equally spaced with width 
$2^{k-\mu}$ and  possible deviation $C$.
Then we assume that there is  $C_1>1$ such that 
$$\sum\Sb J\in \cJ^k_\mu \endSb
 \card (J)\le C_1 N(E^k, 2^{k-\mu})
\tag 1.10$$
for every $k\in \bbZ$, $\mu\in\Bbb N$.

\endroster
\endproclaim


Note that  by Proposition 1.2 the  hypothesis
$(\cR_{p})$
  insures  that the maximal operator associated to the set of endpoints,  $\cup_{k>0} \cD^k$,
 maps $L^p$ to $L^{p}$.


Our main results  are

\proclaim{Theorem  I}
Suppose
 that $1<p<d/(d-1)$
and suppose that $E$
satisfies the
 regularity assumption $(\cR_{p})$.
Then $\cM_E$ is bounded on  $L^p(\bbR^d)$
if and only if  condition $(\cC_{p,p})$ holds.
\endproclaim

\proclaim{Theorem  II}
Suppose
 that $1<p<d/(d-1)$,
and suppose that $E$
satisfies the
 regularity assumption $(\cR_{p})$.
Then $\cM_E$ is of weak type $(p,p)$
if and only if  condition $(\cC_{p,\infty})$ holds.
\endproclaim

It is well known that the weak type $(p,p)$ bounds imply pointwise
convergence theorems. By the Theorems of Calder\'on and Stein
\cite{15 ,ch.X,\S2} and the fact that $(\cC_{p,\infty})$ is necessary for the
$L^p\to L^{p,\infty}$ inequality, these are sharp:

\proclaim{Corollary 1.3} Let $\{t_j\}_{j=1}^\infty$ be a sequence with $\lim_{j\to\infty} t_j=0$
and assume that $E=\{t_j\}$ satisfies condition $(\cR_{p})$ for some
$p\in (1,d/(d-1))$.  Let
$\cA_{t_j} f(x)= \int_{S^{d-1}}f(x-t_jy')d\sigma(y')$.

(i) Suppose that $E$ satisfies condition $(\cC_{p,\infty})$. Then
$\lim_{j\to\infty}\cA_{t_j}f(x)=f(x)$ almost everywhere.

(ii) If  condition $(\cC_{p,\infty})$ is not satisfied then there is a
nonnegative function $f\in L^p(\Bbb R^d)$ such that
$\limsup_{j\to \infty} \cA_{t_j} f(x)=\infty$ almost everywhere.



\endproclaim

\remark{Remark} It would of course be interesting to know whether some 
regularity assumption is needed. As a typical example where the regularity
assumption fails consider the Cantor 
middle third set, translated by $1$, so that
$E_0 =\{1+\sum_{\nu=1}^{\infty} b_{\nu}
3^{-\nu}:  b_{\nu} \in \{0,2\}\}$ and let $E=\cup_{k\in \Bbb Z} E_0$. 
Now the critical exponent is  $p_{\text{cr}}
 = 1 + (d-1)^{-1}
\log 2/\log 3$. The set 
$E_0$ satisfies condition
 ($\cC_{p_{\text{cr}},p_{\text{cr}}}$) and the set $E$ satifies condition
 ($\cC_{p_{\text{cr}},\infty}$). However 
$\cR_{p_{\text{cr}}}$ fails to hold and thus Theorems I and  II above
do not apply.
It is not known whether $\cM_{E_0}$ or $\cM_E$ are of weak type 
$(p_{\text{cr}}, p_{\text{cr}})$ ; see however a counterexample for a closely related maximal 
operator in \S8.2 below.
A much easier result is that 
$\cM_E$ is  of {\it restricted} weak type, see Proposition 1.4 below.
\endremark

We now turn to the limiting case $p=p_d:=d/(d-1)$. There  are sharp results,
at least for $L^{p_d}$ boundedness, although
conditions $(\cR_p)$ and $(\cC_{p,p})$
are replaced by the following different conditions
$(\widetilde\cR_{p_d})$ and 
$(\widetilde\cC_{p_d})$, 
respectively.

\proclaim{Regularity hypothesis $(\widetilde \cR_{p_d})$ }

$E$ satisfies hypothesis $(\widetilde \cR_{p_d})$ if
for each $k$ there is a collection $\cJ^k=\{J\}$  
of subsets of $[2^k,2^{k+1}]$ so that
assumptions (a) and (b) in $(\cR_p)$ hold  but (c) in $(\cR_p)$ is replaced by
\roster
\item"{($\widetilde c$)}"
 There is a $C_1>1$ such that 
$$\sum_{\mu\ge n} 2^{-\mu}\sum\Sb J\in \cJ^k_\mu \endSb
\card (J)\le C_1 2^{-n} N(E^k, 2^{k-n})
\tag 1.11$$
holds uniformly in $n\in\Bbb N$.
\endroster
\endproclaim

The analogue of condition $(\cC_{p,p})$
is
\proclaim{Condition $(\widetilde \cC_{p_d})$} 


The discrete measure $\sum_{k\in \Bbb Z}
\sum_{n>0} N(E^k, 2^{k-n}) 2^{-n} n^{1/(d-1)}\delta_{k,n}$
 is a Carleson measure on the upper half plane; i.e.
$$
\sup_{|I|\ge 1}\frac 1{|I|}\sum_{(k,n)\in T(I)}
N(E^k,2^{k-n})2^{-n} n^{1/(d-1)}
<\infty
\tag 1.12
$$
where the supremum is taken over all intervals of length $\ge 1$ and
$T(I)$ is the tent of $I$, i.e.
$T(I)=\{(x,t):x\in I,\, 0\le t\le|I|\}.$
\endproclaim




It was shown in \cite{13} that for $d\ge 3$  condition 
$(\widetilde \cC_{p_d})$ is equivalent
with the $L^{p_d}$ boundedness of $\cM_E$ on {\it radial} functions.
For general $L^p$ functions we have a  similar result 
provided that hypothesis $(\widetilde \cR_{p_d})$ is satisfied:



\proclaim{Theorem III}
Let $d\ge 2$ and  $p_d=d/(d-1)$ and suppose
  that $E$
satisfies the
 regularity assumption $(\widetilde \cR_{p_d})$.
Then $\cM_E$ is bounded on  $L^{p_d}(\bbR^d)$
if and only if  condition $(\widetilde \cC_{p_d})$ holds.
\endproclaim

Concerning a weak type $(p_d,p_d)$ inequality in dimensions $d\ge 3$ one may conjecture
that the hypothesis
$$
N(E^k, 2^k\delta)\le C\delta^{-1} \bigl[\log(1/\delta)\bigr]^{-1/(d-1)}
\tag 1.13
$$
is necessary and sufficient for $L^{p_d}\to L^{p_d,\infty}$ boundedness as this is shown to hold in \cite{13} on $L^{p_d}_\rad$.
For general functions $f$ and under the regularity assumption $(\widetilde \cR_{p_d})$
we prove the following slightly weaker result.

\proclaim{Theorem IV}
Let $d\ge 3$ and $p_d=d/(d-1)$ and suppose
  that $E$
satisfies the
 regularity assumption $(\widetilde \cR_{p_d})$ and suppose that
$$
N(E^k, 2^k\delta)\le C\delta^{-1} \bigl[\log(1/\delta)\bigr]^{-1/(d-1)}
\bigl[\log\log (1/\delta)\bigr]^{-1}
\tag 1.14
$$
uniformly in $k\in \bbZ$ and $\delta\le e^{-2}$. 
Then $\cM_E$ is of weak type $(p_d,p_d)$.
\endproclaim


At present we do not know whether  the same conclusion holds 
under the weaker condition
(1.13). This accounts for
the as yet undecided weak type $(p_d,p_d)$ estimate for
$M_{\widetilde E(\beta)}$ in the remaining case $\beta = 1/(d-1)$ in Theorem 1.1.


We now briefly  turn to the question of {\it restricted weak type inequalities}.
Here no regularity assumption is needed.

\proclaim{Proposition 1.4}
Let $1<p\le d/(d-1)$, $d\ge 3$ or $1<p<2$, $d=2$ and suppose that  
$E$ satisfies 
condition ($\cC_{p,\infty}$). Then $\cM_E$ is of restricted 
weak type $(p,p)$, i.e.
it maps $L^{p,1}$ to $L^{p,\infty}$.
\endproclaim

It remains open whether for the range $1<p<d/(d-1)$ the operator is of weak type 
$(p,p)$, under condition
($\cC_{p,\infty}$) alone,  without the regularity assumption.
Proposition 1.4  is much more straightforward  than Theorem II above and we shall not give the details of the 
proof here. 
 For $p_d=d/(d-1)$, $d\ge 3$  the result  had been
already proved by Bourgain  \cite{1}, and a variant of his argument applies  for $1<p<d/(d-1)$
 as well.
Indeed  let $\cA_t^j$ be the frequency localized  operator as in  (2.1) 
below and define
the maximal operator $M_j$ by  $M_jf(x)=\sup_{t\in E} |\cA^j_t f(x)|$. Then the estimates in
\cite{12} show that for $1<q\le 2$ the operator $M_j$ is bounded on $L^q$ with norm
$O(2^{j(d-1)(\frac 1{q'}-\frac{p-1}q)})$ and the argument in \cite{1} shows the
restricted weak type estimate. The argument fails for $p=d=2$ and in fact the question
whether the full circular maximal function  is of restricted weak type (2,2)
({\it i.e.} maps $L^{2,1}$ to $L^{2,\infty}$) had been posed in
\cite{16}. We note that Leckband \cite{7} proved  that for
{\it radial} functions one has indeed $L^{2,1}_\rad\to L^{2,\infty}$ 
boundedness.
However a Besicovitch set construction
can be used to disprove
 the restricted weak type (2,2) inequality for general functions.
The argument (see \S8 below)  shows

\proclaim{Proposition 1.5}
Suppose $d=2$ and
$$\sup_{k>0} \sup_{\delta<1/10}
N(E^k, 2^k\delta) \delta \log\delta^{-1}
=\infty.
$$
Then
$\cM_E$ is not of restricted weak type $(2,2)$.
\endproclaim

{\it Structure of the paper:}
In \S2 we shall review
some essentially known estimates for spherical means which are needed later.
In \S3 we shall review atomic decompositions in $L^p$.
\S4 contains a proof of the $L^p$ estimates as stated in Proposition  1.2 and Theorem I.
The weak type $(p,p)$  inequalities (Theorem II) are proved in \S5.
The necessary modifications for the proofs of Theorem III and IV 
are discussed in 
\S6 and \S7, respectively.
In \S8 we  discuss some examples and include  the proof of Proposition 1.5.


\head{\bf 2. Estimates on spherical means}
\endhead



We shall need to introduce regularizations of $\cA_t$ in (1.1) via
dyadic frequency  cutoffs.
Let $\beta_0$ be a radial $C^\infty_0$ function so that
$\beta_0(\xi)=1$ if $|\xi|\le 1$ and $\beta_0(\xi)=0$ if $|\xi|\ge 2$.
For $j=1,2,\dots$ let $\beta_j(\xi)=\beta_0(2^{-j}\xi)-\beta_0(2^{1-j}\xi)$
and define
$\cA^j_t$ by
$$ \widehat {\cA^j_t f}(\xi)= \widehat{d\sigma}(t\xi)
\beta_j(t\xi)\widehat f(\xi)
\tag 2.1
$$
so that
$$\cA_t=\sum_{j=0}^\infty \cA^j_t.
$$

Let $\tilde \beta$ be a radial $C^\infty_0$ function  which is supported
where $2^{-6}\le |\xi|\le 2^{6}$
 and equal to $1$ when
$2^{-5}\le |\xi|\le 2^{5}$. Let $P^l\! f$ be defined by
$\widehat{P^l\! f} (\xi)=\tilde \beta(2^{-l}\xi)
\widehat f(\xi)$ and observe that
$$
\cA^j_t f= \cA^j_t P^{j-k} f \quad\text{ if } t\in E^k.
\tag 2.2
$$

Clearly the maximal function $\sup_{t>0} |\cA^j_t f(x)|$ is dominated by
$C_j M_{HL} f(x)$ where $M_{HL}$ is
the Hardy-Littlewood maximal function of $f$;
in fact $C_j=O(2^j)$ (\cf. Lemma 2.1 below).
Therefore
$$
\cM_E f(x)\le  M_{HL} f(x)+ \sup_{k\in\bbZ}\sup_{t\in E^k }
|\sum_{j\ge 10}  \cA^j_t  P_{j-k} f(x)|,
\tag 2.3
$$
and throughout this paper we shall assume that summations in $j$ are
extended over  $j\ge 10$.

Here we collect well known estimates on spherical means and its regularization
$\cA^j_t$  which were  used in this
or a related form in previous papers (in particular see \cite{12} for some of the more  technical
statements).

\proclaim{Lemma 2.1} Let $\cA^j_t$ be as above and let $B_t^j=\frac {d}{dt}\cA_t^j$.
Suppose that  $2^k\le t\le 2^{k+1}$, $j\ge 10$  and that $1\le p\le 2$.

(i)
$$|\cA^j_t f(x)|
+2^{-j}t|B_t^j f(x)|
\le C_M 2^j \int \frac{t^{-d}}{(1+2^j|\frac{|x-y|}{t} -1|)^M} |f(y)|dy
$$

(ii)
$$\|\cA^j_t\|_{L^p\to L^p}+
2^{k-j}\|B^j_t\|_{L^p\to L^p}
\lc 2^{-j(d-1)/p'}.$$


(iii) Let $I\subset [2^{k-1}, 2^{k+2}]$ be an interval of length
$2^{k-j}$. Then
$$
\big\| \sup_{t\in I}|\cA^j_t f|\big\|_{L^p}
\lc
 2^{-j(d-1)/p'}\|f\|_{L^p}.
$$


\endproclaim

\demo{Sketch of Proof}  (i) is a straightforward  calculation, which also implies (ii) for $p=1$.
It is well known that
$|\widehat {d\sigma}(\xi)|\lc (1+|\xi|)^{-(d-1)/2}$ and
thus (ii) for $p=2$ follows, and interpolation settles the case $1<p<2$.
(iii) follows by writing  $\cA_t^j=\cA_{t_0}^j+\int_{t_0}^t B_s^j ds$ for $t_0\in I$.
\enddemo

\definition{Definition}
For a set $\cE$ of dilations and $L\in \bbZ$, let
 $\fI_{L}(\cE)$ be a minimal collection of dyadic  intervals
of length $2^{L}$ covering $\cE$. For $I\in \fI_L(\cE)$ let $r_I$
denote the midpoint of the
interval $I$, and  for a dyadic cube $Q$, let  $2^{L(Q)}$ denote its
sidelength. Then for $\eta \ge 1$, we define
$$V_{Q,\eta}(\cE)=\bigcup_{I\in\fI_{L(Q)}(\cE)}
\{x\in \Bbb R^d: \big||x-x_Q|-r_I\big|\le 2^{L(Q)+4}\eta\};
\tag 2.4
$$
for $\eta=1$ we also write
$V_{Q}(\cE)=
V_{Q,1}(\cE)$.
\enddefinition

\proclaim{Lemma 2.2} Let $\cE \subset [2^k,2^{k+1}]$.

(i)
For  $1\le p\le 2$,
$$
\big\|\sup_{t\in \cE} |\cA^j_t f|\big\|_{L^p}\lc [N(\cE,2^{k-j})]^{1/p} 
2^{-j(d-1)/p'}\|f\|_{L^p}.
$$

(ii)
Let $Q$ be a dyadic cube, let
 $f_Q$ be an $L^2$ function supported on  $Q$ and suppose
$k-j \le L(Q) \le k-10$.
Then
 $$\big\|\sup_{t\in \cE}|\cA^j_t f_Q|\big\|_{L^1(V_{Q}(\cE))}\lc
2^{(-L(Q)+k-j)(d-1)/2}
N(\cE,2^{k-j})
2^{L(Q)d/2}\|f_Q\|_{L^2} .
$$

(iii) Let
$\cQ$ be a collection of pairwise disjoint cubes of
sidelength $2^{k-j+\sigma }$ where $\sigma \ge 0$.
Then for $\sigma\le j$,
 $$\Big\|\sup_{t\in \cE}\big|\cA^j_t [\sum_{Q\in\cQ} f_Q]\big|\Big\|_{L^p}
\lc
2^{-\sigma (d-1)(1/p-1/2)}
[N(\cE,2^{k-j})]^{1/p} 2^{-j(d-1)/p'}
\Big(\sum_{Q\in\cQ}
|Q|^{1-p/2}\|f_Q\|_{L^2}^p\Big)^{1/p} .
$$

(iv)
Let $\cQ$ be as in (iii) and let $\cV$ be an open set 
containing $\bigcup_{Q\in \cQ} V_{Q,\eta}(\cE)$.
Then, for $\eta\ge 1$,
$$\Big\|\sup_{t\in \cE}\big|\cA^j_t[\sum_{Q\in \cQ} f_Q]\big|\Big
\|_{L^p(\Bbb R^d \setminus \cV)}\le C_M
(2^{\sigma  }\eta)^{-M(\frac 2p-1)}
[N(\cE,2^{k-j})]^{1/p} 2^{-j(d-1)/p'}
\Big(\sum_{Q\in \cQ}\|f_Q\|_{L^p}^p\Big)^{1/p} .
$$


(v) The estimates in (i), (ii), (iii) and (iv) remain valid if for $t\in \cE$
the operator $\cA^j_t$ is replaced by $2^{k-j}B^j_t=
2^{k-j}\frac{d}{dt}\cA^j_t$.
\endproclaim

\demo{Sketch of Proof}
(i) is a rather straightforward consequence of Lemma 2.1, (iii).
To prove (ii) we use Cauchy-Schwarz
to pass from an $L^1$ estimate
on the exceptional set $V_Q({\Cal E})$ to an $L^2$ estimate
 (namely Lemma 2.1 (ii) with $p=2$), 
and for the estimate
off the exceptional set we use the explicit form (2.4).
(iv) for $p=2$ is a consequence of (i), and (iv) 
for $p=1$ follows from the explicit form of the kernel in Lemma 2.1 (i).
The general case is obtained by
 interpolation. 
(iii) for $p=2$ is a consequence of (i), and (iii) for $p=1$
follows from (ii) and (iv). The general case is obtained by interpolation.
\qed
\enddemo

A small variant is

\proclaim{Lemma 2.3}
Let $J\subset [2^k,2^{k+1}]$ be an equally spaced set
with width $2^{k-\mu}$ (here $\mu\ge 0$)  and possible deviation $B$, and let $a_J<b_J$
be the endpoints of $J$. Suppose that $b_J-a_J\ge 2^{k-j}$ and
$\mu\ge j$.
Then the following statements hold.

(i) $$
\Big\|\sup_{t\in J}|\cA^j_t f|\Big\|_{L^p}
\le C_B N(J, 2^{k-\mu})^{1/p} 2^{-j(d-1)/p'}2^{(j-\mu)/p} \|f\|_{L^p}.
$$

(ii) Let
$\cQ$ be a collection of pairwise disjoint cubes of
sidelength $2^{k-j+\sigma }$ where $\sigma \ge 0$.
Then for $\sigma \le j$,
 $$
\Big\|\sup_{t\in J}\big|\cA^j_t [\sum_{Q\in\cQ} f_Q]\big|\Big\|_{L^p}
\le C_B
2^{-\sigma (d-1)(1/p-1/2)}
\card(J)^{1/p}
 2^{-j(d-1)/p'} 2^{(j-\mu)/p}
\Big(\sum_{Q\in\cQ }
|Q|^{1-p/2}\|f_Q\|_{L^2}^p\Big)^{1/p}.
$$

(iii) Let $\cQ$ be as in (ii) and let
 $\cV$ be an open set containing $\bigcup_{Q\in \cQ} V_{Q,\eta}(\cE)$.
Then, for $\eta\ge 1$,
$$
\Big\|\sup_{t\in J}\big|\cA^j_t
[\sum_{Q\in \cQ} f_Q]\big|\Big
\|_{L^p(\Bbb R^d \setminus \cV)}
\le C_{B,M}
(2^{\sigma  }\eta)^{-M(\frac 2p-1)}
\card(J)^{1/p}
 2^{(j-\mu)/p} 2^{-j(d-1)/p'}
\Big(\sum_{Q\in \cQ}\|f_Q\|_{L^p}^p\Big)^{1/p}.
$$

\endproclaim

\demo{Proof} We simply observe that if $b_J-a_J\ge 2^{k-j}$
then 
$$N(J, 2^{k-j})\approx 2^{j-\mu} N(J, 2^{k-\mu})\approx 2^{j-\mu}\card(J)
\tag 2.5
$$
and the conclusions (i)-(iii) follow from Lemma 2.2.
\qed \enddemo


\head{\bf 3. Atomic decompositions
}
\endhead

We give a decomposition of the maximal operator and also the 
function it acts on;
this is motivated by one of the proofs of the standard 
atomic decomposition (following \cite{3}, \cite{9}) based on 
square functions;
used for example in the theory of Hardy spaces on product domains.

 For $c_0=10\sqrt d$ let
$$\cN^k f(x)=\sup_{|y|\le c_0 2^{-k}} |P^kf(x+y)|$$ and define the maximal
square function
$$\cN f(x)=\Big(\sum_{k=-\infty}^\infty|\cN^k f(x)|^2\Big)^{1/2};$$
then  $$\|\cN f\|_{L^p}\approx_p \|f\|_{L^p}, \qquad 1<p<\infty,
\tag 3.1$$  and
 $\|\cN f\|_{L^1}\approx \|f\|_{H^1}$.

Consider the level sets $\Omega_n=\{x: \cN f(x)>2^n\}$ and the expanded sets
$\widetilde \Omega_n=\{ x: M_{HL}\chi_\Omega(x)>1/2\}$; 
here $M_{HL}$ is the Hardy-Littlewood
maximal function. Then $|\widetilde \Omega_n|\le C|\Omega_n|$. 
Let $\fR$ denote the family
of all dyadic cubes
and let $\fR_n$, for $n\in \Bbb Z$, denote the collection of all dyadic
cubes $R$  with the property that $|R\cap\Omega_n|>|R|/2$ but
$|R\cap\Omega_{n+1}|\le|R|/2$. Then from these definitions one easily deduces
$$
\sum_{k=-\infty}^\infty\sum\Sb R\in \fR_n\\L(R)=-k\endSb
\|(P^k f)\chi_R\|_{L^2}^2\lc 2^{2n}|\Omega_n|
\tag 3.2
$$
 (see for example  Lemma 3.1 in \cite{9}).

Let $e_R= (P^l f)\chi_R$ if $L(R)=-l$. Then from (2.2) we have 
$$ 
\cA_t^j f = \cA_t^j P^{j-k} f = 
\cA_t^j[\sum_{L(R)=k-j}
e_R]
\tag 3.3
$$
if $t\in E^k$.

Now  $\widetilde \Omega_n$ is an open set with finite 
measure and we can form the 
Whitney-decomposition into dyadic cubes. 
Let $\fW_n$ be the set of Whitney cubes and observe that every
$R\in \fR_n$ is contained in a unique Whitney-cube $Q(R)$. 
This defines a function $R\mapsto Q(R)$ for all dyadic cubes.

For a dyadic cube $Q$ we define now
$$F_Q^l(f)=\sum\Sb Q(R)=Q\\L(R)=-l\endSb e_R.
\tag 3.4
$$
Notice that
$F_Q^l=0$ if $-l>L(Q)$.

>From (2.3) and (3.3), we have the pointwise estimate 
$${\cM}_E f(x) \lc M_{HL} f(x) + 
\sup_{k\in \Bbb Z}\sup_{t\in E^k} \sum_{j\ge 10}\Big| \cA_t^j[\sum_{L(Q)\ge 
k-j}
 F_Q^{j-k}(f)](x)\Big|.
\tag 3.5
$$


It is useful to introduce a space $X^p$ of vector-valued functions
as follows.
\definition{Definition}
Let $X^p$ be the space of
vector-valued functions $F=(F^l_Q)$ where the dyadic cubes $Q$ satisfy
$L(Q)+l\ge 0$,
$F^l_Q$ is supported  on $Q$,
and
$$
\|F\|_{X^p}=\Big(\sum_Q|Q|^{1-p/2}\big(\sum_{l:L(Q)+l\ge 0}\|F^l_Q\|_{L^2}^2
\big)^{p/2}\Big)^{1/p}
\tag 3.6
$$
is finite.
\enddefinition

We first observe

\proclaim{Lemma  3.1} For $1\le p\le 2$, 
$$\|F(f)\|_{X^p}\lc \|\cN f\|_{L^p}.$$
\endproclaim

\demo{Proof}
We write
$$
\align\|F(f)\|_{X^p}&=\Big(\sum_Q|Q|^{1-p/2}\Big(
\sum_{\ell: L(Q)+\ell\ge 0} \Big\|
\sum\Sb Q(R)=Q\\L(R)=-\ell \endSb e_R
\Big\|_{L^2}^2\Big)^{p/2}\Big)^{1/p}
\\
&\le \Big(\sum_Q|Q|^{1-p/2}\Big(\sum_{\ell: L(Q)+\ell\ge 0}
\sum_n
\sum\Sb R\in \fR_n\\ Q(R)=Q\\L(R)=-\ell \endSb
\|e_R\|_{L^2}^2\Big)^{p/2}\Big)^{1/p} .
\endalign
$$
Now we use the imbedding $\ell^p\subset\ell^2$ for $p\le 2$
 to estimate the last expression by
$$
\align
&\Big(\sum_Q|Q|^{1-p/2} \sum_n\Big(\sum_{\ell: L(Q)+\ell\ge 0}
\sum\Sb R\in \fR_n\\ Q(R)=Q\\L(R)=-\ell \endSb
\|e_R\|_{L^2}^2\Big)^{p/2}\Big)^{1/p}
\\
&\le
 \Big(\sum_n\sum_{Q\in \fW_n}|Q|^{1-p/2} \Big(
\sum\Sb R\in \fR_n\\ Q(R)=Q \endSb \|e_R\|_{L^2}^2\Big)^{p/2}\Big)^{1/p}
\endalign
$$
and
by (3.2) and several
applications of H\"older's inequality
this in turn is estimated by
$$
\align
&\Big(\sum_n \big(\sum_{Q\in \fW_n}|Q|\big)^{1-p/2} \Big(\sum_Q
\sum\Sb R\in \fR_n\\ Q(R)=Q \endSb \|e_R\|_{L^2}^2\Big)^{p/2}\Big)^{1/p}
\\
&\le
\Big(\sum_n |\widetilde \Omega_n|^{1-p/2} \Big(\sum\Sb R\in \fR_n \endSb
\|e_R\|_{L^2}^2\Big)^{p/2}
\Big)^{1/p}
\\
&\lc
\Big(\sum_n |\widetilde \Omega_n|^{1-p/2} \big(2^{2n}|\Omega_n|\big)^{p/2}
\Big)^{1/p}
\\
&\lc
\Big(\sum_n | \Omega_n|2^{np}\Big)^{1/p}  \lc \|\cN f\|_{L^p} .
\endalign
$$
This proves the Lemma.\qed
\enddemo

We now return to estimate  the second term on the right of (3.5).
The part where the sum extends over
cubes $Q$ with $L(Q) \le k$ is the most difficult to handle.
In the following lemma we  shall first dispose of the remaining part which is dealt with 
by straightforward $L^2$ estimates.


\proclaim{Lemma 3.2}
Let $1<p\le 2$, $\eps>0$ and suppose that 
$$\sup_k N(E^k,2^k\delta)^{1/2} \delta^{(d-1-\eps)/2}\le A.$$
Let  $\{\chi_{Q,l}\}$ be a family of measurable functions so that
$$\sup_l\|\sum_Q|\chi_{Q,l}\|_{L^2} \le 1\tag 3.7$$ and define
$$
\fN_j F(x)=
\sup_{k\in \Bbb Z}\sup_{t\in E^k} \Big| \cA_t^j[\sum_{L(Q)\ge
k} \chi_{Q,k}
 F_Q^{j-k}](x)\Big|.
\tag 3.8
$$
Then
$$ \| \fN_j F \|_{L^p} \le C 2^{-\eps j}  A \|F\|_{X^p} \tag 3.9 $$
where $C$ is independent of the choice of the particular family $\{\chi_{Q,k}\}$
.
\endproclaim

\demo{\bf Proof} We shall verify (3.9) for $p=1$ and for $p=2$; the general case follows by 
interpolation.

For $p=2$ we replace the $\sup$ in $k$ by a square function and  use Lemma 2.2 (i) to obtain
$$
\align
\|\fN_j F\|_{L^2}&\le
\Big(\sum_k\Big\|\sup_{t\in E^k}| \cA_t^j \big[ \sum\Sb Q:\\L(Q)\ge k \endSb
\chi_{Q,k}
F_Q^{j-k}
\big]\Big\|_{L^2}^2\Big)^{1/2}
\\
&\le C A 2^{-\eps j}
\Big(\sum_k\Big\| \sum\Sb Q:\\L(Q)\ge k \endSb \chi_{Q,k} F_Q^{j-k}
\Big\|_{L^2}^2\Big)^{1/2}
\\
&\le 
C A 2^{-\eps j}
\Big(\sum_k\sum\Sb Q:\\L(Q)\ge k \endSb \Big\|F_Q^{j-k}
\Big\|_{L^2}^2\Big)^{1/2}
\endalign
$$
where for the last inequality we have used the assumption on the family $\{\chi_
{Q,l}\}$.
This proves (3.9) for $p=2$.

Now consider the case $p=1$. Given a cube $Q$ we let $Q^*$ denote the cube with same center but tenfold sidelength.
We then estimate (following standard procedure in estimations of singular integrals acting 
on atoms)
$$
\|\fN_j F\|_{L^1}\le  \sum_Q I_Q+II_Q
$$
where
$$\align
I_Q&=
\big\|\sup_k\sup_{t\in E^k}| \cA_t^j [\chi_{Q,k} F_Q^{j-k}]
\big\|_{L^1(Q^*)}
\\
II_Q &=
\big\|\sup_k\sup_{t\in E^k}| \cA_t^j [\chi_{Q,k} F_Q^{j-k}]
\big\|_{L^1(\bbR^d\setminus Q^*)}.
\endalign
$$

Now for $I_Q$ we use the Cauchy-Schwarz inequality and the $L^2$ estimate above to deduce
that
$$\align I_Q&\lc |Q|^{1/2}
\big\|\sup_k\sup_{t\in E^k}| \cA_t^j [\chi_{Q,k} F_Q^{j-k}]
\big\|_{L^2}
\\
&\le C A 2^{-\eps j} |Q|^{1/2}
\Big(\sum_k\big\|F_Q^{j-k}
\big\|_{L^2}^2\Big)^{1/2}
\endalign
$$

For $II_Q$ we use Lemma 2.1 (i). 
In that formula we use that if $y\in Q$, $x\in Q^*$, $L(Q)\ge k$, $t\le 2^{k+1}$ then
$|t^{-1}|x-y|-1|\approx t^{-1}|x-y|$ and thus for $M>d\ge 2$
$$
\align
II_Q&\lc 2^{j} \sum_{k\le L(Q)}\int_{|x-y_Q|\ge 2^{L(Q)+2}} 2^{k(M-d)} 2^{-jM}\int_Q|x-y|^{-M} 
|F^{j-k}_Q(y)| dy dx
\\ &\lc 2^{j(1-M)}\sum_{k\le L(Q)} 2^{(k-L(Q))(M-d)} \|F^{j-k}_Q\|_{L^1(Q)}
\\
&\lc 2^{j(1-M)} |Q|^{1/2} \Big(\sum_k \|F^{j-k}_Q\|_{L^2}^2\Big)^{1/2}.
\endalign
$$
Now  $M$ can be chosen to be $\ge 1+\eps$ and we obtain 
that
$\sum_Q (I_Q+II_Q)$ is bounded by $C A 2^{-\eps j} \|F\|_{X^1}$, 
thus proving (3.9) for $p=1$.
\qed\enddemo

For the remainder of the paper we will only have to deal with the part in (3.5) dealing with the contribution $k>L(Q)$.
 Define for a positive  integer $\sigma$
$$\fM_\s F(x)= \sup_k\sup_{t\in E^k}
\Big|\sum_{j\ge \sigma}  \cA^j_t [\sum\Sb Q:\\ 
L(Q)= k-j+\s\endSb F^{j-k}_Q ](x)\Big
|.
\tag 3.10
$$

Our main reduction in this section is

\proclaim{Proposition 3.3}
Let  $1<p< 2$, suppose that hypothesis ($\cC_{p,\infty}$) is satisfied  and 
suppose that for some $\eps_0>0$ the inequality
$$\|\fM_\s F\|_{L^{p,q}}\le C_0 2^{-\eps_0 \s}\|F\|_{X^p}
\tag 3.10
$$
holds for all compactly supported $F$
(meaning that $F^l_Q$ vanishes for all but finitely many
$l$ and $Q$). Then   there is $c(p,\eps_0)>0$ so that
$$\|\cM_E f\|_{L^{p,q}} \le c(p,\eps_0)C_0 \|f\|_{L^p}
$$
for all $f\in L^p(\bbR^d)$.
\endproclaim

\demo{\bf Proof}
Let $ F_Q^l(f)$ be as in (3.4).
For $\sigma=1,2,,\dots$ define
$F^{(1)}_\s(f)$ by
$[F^{(1)}_\s]_Q^l(f)= F_Q^l(f)$ if
$L(Q)=\s-l$ and
$[F^{(1)}_\s]_Q^l(f)=0$
if $L(Q)\neq \s-l$.
For $j\ge 10$ define
$[F^{(2)}_j]_Q^l(f)= F_Q^{j+l}(f)$ if $L(Q)\ge -l$ and zero otherwise
and let
$\chi_{Q,l}^j$ be the characteristic function of
$\cup_{n\in \bbZ}\cup\{R: R\in \fR_n, L(R)=-l-j, Q(R)=Q\}$.
Then for every fixed  $j$ condition (3.7) is satisfied for the family  
$\{\chi_{
Q,l}^j\}$.

>From  (3.5) we get
$${\cM}_E f(x) \lc M_{HL} f(x) +  \sum_{\sigma>0}
\fM_\sigma[ F_\s^{(1)}(f)](x) 
+  \sum_{j\ge 10}
\fN_j[ F_j^{(2)}(f)](x) .
$$
Note that it follows from Lemma 3.1 and (3.1) that
$\|F^{(1)}_\s(f)\|_{X^p}\le C_p \|f\|_p$
and 
$\|F^{(2)}_j(f)\|_{X^p}\le C_p \|f\|_p$ for
$1<p\le 2$, uniformly in $\sigma$ and $j$. 
>From hypothesis ($\cC_{p,\infty}$) it follows that the assumption of Lemma 3.2 holds with 
$\eps=(d-1)(2-p)$ which is positive since we are assuming $p<2$.  Thus 
$$\Big \|
\sum_{j\ge 10}
\fN_j[ F^{(2)}_j(f)]\Big\|_{L^{p,q}}\lc
\sum_{j\ge 10} 2^{-\eps j}
\| F^{(2)}_j(f)\|_{X^p}\lc
 \|f\|_p.
$$
By our assumption we also have
$$\Big \|
\sum_{\s>0}
\fM_\sigma[ F^{(1)}_\s(f)]\Big\|_{L^{p,q}}
\lc\sum_{\s>0} 2^{-\eps_0\sigma}\| F^{(1)}_\s(f)\|_{X^p}\lc
 \|f\|_p
$$
and the proposition is proved.
\qed
\enddemo

\head{\bf 4. $\boldkey L^{\boldkey p}$ estimates}\endhead
We shall use Proposition 3.3 and in order
to  prove  $L^p$ estimates we have to verify the $X^p\to L^p$ estimate
for $\fM_\s$ in (3.10). 
We shall first prove Proposition 1.2 where no regularity assumption is needed.

We shall also use the following definitions.
$$
\align
G_\sigma ^l (F)& =
\sum_{Q: L(Q)=-l+\sigma } F^{l}_Q
\tag 4.1
\endalign
$$
and let
$G_\sigma (F)=\{G^l_\sigma  (F)\}_{\l\in\bbZ}$ 
be the corresponding vector valued analogue.
\proclaim{Proposition 4.1} Suppose that  $1< p<2$ and suppose that
$\sum_n\omega_n^{-p'}\le 1$. Let
$$
\la_{j,k}= N(E^k, 2^{k-j})^{1/p} 2^{-j(d-1)/p'} .
\tag 4.2
$$
Then 
$$
\|\fM_\sigma F\|_{L^p}\lc 
 2^{-\sigma(d-1)(1/p-1/2)}
\sup_l \Big(\sum_{n=0}^\infty|\om_n|^p \la_{n,l+n}^p\Big)^{1/p}
\|F\|_{X^p} .
$$\endproclaim
\demo{\bf Proof}
We  estimate
using H\"older's inequality
$$\align
|\fM_\sigma  F (x)|&\lc
\Big(\sum_k\big \|\{\om_{\cdot}^{-1}\}\big\|_{\ell^{p'}}^p  \sum_{j=10}^\infty
\omega_j^p
\sup_{t\in E^k}|\cA^j_{t}  G^{j-k}_\sigma ( F)|^p\Big)^{1/p}.
\endalign
$$
By Lemma 2.2 (iii) the $L^p$ norm of this expression is estimated as
$$
\align
\|\fM_\sigma  F\|_{L^p} &\lc 2^{-\sigma (d-1)(1/p-1/2)}
\Big(\sum_k\sum_j\big[
 \omega_j \la_{j,k}\big]^p \sum\Sb L(Q)=k-j+\sigma \endSb
|Q|^{1-p/2}\big\|
 F^{j-k}_Q
\big \|_{L^2}^p\Big)^{1/p}
\\&\lc 2^{-\sigma (d-1)(1/p-1/2)}
\Big(\sum_Q\sum_j \big[\om_j \la_{j,L(Q)+j-\sigma }
 \big]^p
|Q|^{1-p/2}\big\|
F^{\sigma-L(Q) }_Q
\big \|_{L^2}^p\Big)^{1/p}
\\
&\lc 2^{-\sigma (d-1)(1/p-1/2)}
\Big(\sum_l \sum_{Q: L(Q)=\sigma-l }
\sum_j [\om_j  \la_{j, l+j}]^p
|Q|^{1-p/2}\big\|
F^{l}_Q
\big \|_{L^2}^p\Big)^{1/p}
\\
&\lc 2^{-\sigma (d-1)(1/p-1/2)} \|F\|_{X^p}. \ \qed
\endalign
$$
\enddemo

\demo{\bf Proof of Proposition 1.2} Immediate from Propositions 3.3 and 4.1
when $1<p<2$. The case $p=2$ (and hence $d=2$) follows as in the
proof of Proposition 4.1 where now we treat the whole operator
$\cM_E$. \enddemo

We now turn to the proof of  the
\subheading{$\boldkey L^{\boldkey p}$ estimates under  the
regularity hypothesis}
For the remainder of this section we 
%
shall fix  a choice 
of $\cJ^k$, $\cJ^k_\mu$  as in the definition of
regularity assumption $(\cR_p)$.

Let $\s$ be a positive  integer. Let
$$R_\s F (x)=\sup_k \sup_{J\in \cJ^k}
\sup_{t\in J}
\Big|\sum\Sb j\ge\s :\\
b_J-a_J\le 2^{k-j}\endSb\cA^j_t G^{j-k}_{\s}(F)(x)\Big|
\tag 4.3 $$
and, for 
$m\ge0$
$$
S_{m,\s} F(x)=\sup_{k}\sup\Sb \mu:\\ \mu+m\ge\s\endSb \sup\Sb
J\in\cJ^k_\mu\\ b_J-a_J> 2^{k-\mu-m}\endSb\sup_{t\in  J}
\Big|\cA_t^{\mu+m}
G^{\mu+m-k}_{\s}(F)(x)\Big|.
\tag 4.4$$

Next let
$$M_{\sigma } F(x)\,=\,
\sup_k
\sup_{\mu>\s}\sup_{J\in\cJ^{k}_\mu}\sup_{t\in J}
\big|\sum\Sb \s <j<\mu\\
b_J-a_J> 2^{k-j}\endSb 
\cA_t^{j}G_\sigma ^{j-k}(F)(x)
\big|.
\tag 4.5
$$
Thus
$$
\fM_{\s} F(x)\le R_\s F(x) + \sum_{m} S_{m,\s} F(x) + M_\s F(x).
\tag 4.6
$$

Finally, for $\ell>0$, let
$$\cJ^{k,\ell}_\mu=\{ J\in \cJ^{k}_\mu: b_J-a_J\ge 2^{k-\mu+\ell}\}
\tag 4.7 $$
and define 
$$M_{\ell,\sigma } F(x)\,=\,
\sup_k
\sup_{\mu>\ell + \s}\sup_{J\in\cJ^{k,\ell}_\mu}\sup_{t\in J}
\big|\cA_t^{\mu-\ell}G_\sigma ^{\mu-\ell-k}(F)(x)
\big|
\tag 4.8
$$
so that
$$
M_\s F(x)\le \sum_{\ell>0} M_{\ell,\s} F(x).
\tag 4.9
$$

The claim in Theorem I  will be a consequence 
of the following Propositions 4.2, 4.3, 4.4, 
in conjunction with Proposition 3.3.

The following result is essentially Proposition 4.1 applied to the set of `endpoints', {\it i.e.}
$\cup_k \cD^k$.

\proclaim{Proposition 4.2} 
Suppose that $1\le p\le 2$ and $E$ satisfies
the regularity hypothesis $(\cR_p)$, and let $\cD^k$ be as in (1.9).
Assume that $\{\om_n^{-1}\} \in \ell^{p'}$ with norm $\le 1$. Then
$$\|R_\s F\|_{L^p}\lc 2^{-\s(d-1)(1/p-1/2)}
\sup_l \Big(\sum_{j\ge 0} [\om_j N(\cD^{j+l}, 2^{l})^{1/p}
2^{-j(d-1)/p'}]^p\Big)^{1/p}\|F\|_{X^p} .
\tag 4.10
$$
\endproclaim

\demo{\bf Proof} 
Using H\"older's inequality as above (with $\{\omega_n^{-1}\}\in \ell^{p'})$ we may estimate
$$
R_\s F(x)\lc\Big(\sum_{k}\sum_{j\ge \s}\omega_j^p 
 \sup\Sb J:\\ b_J-a_J\le  2^{k-j}\endSb\sup_{t\in J}|\cA^j_t 
G_\s^{j-k}F(x)|^p\Big)^{1/p} .
$$
Now if for fixed $j,k$ we let
$\cE=\cup_{J\in \cJ^k : b_J-a_J\le 2^{k-j}} J$ then 
$N(\cE, 2^{k-j}) \lc N(\cD^k, 2^{k-j})$. Hence by Lemma 2.2 (iii),
$$\align
\|R_\s F\|_{L^p}
&\lc
\Big(\sum_{k}\sum_{j}\omega_j^p \Big\|
 \sup\Sb J\in \cJ^k :\\ b_J-a_J\le  2^{k-j}\endSb
\sup_{t\in J}|\cA^j_t G_\s^{j-k}(F)|\Big\|_{L^p}^p\Big)^{1/p}
\\
&\lc 2^{-\s(d-1)(1/p - 1/2)}
\Big(\sum_{k}\sum_{j}\omega_j^p N(\cD^k, 2^{k-j}) 2^{-j(d-1)p/p'} 
\sum\Sb L(Q)=k-j+\sigma \endSb
|Q|^{1-p/2}\big\|
 F^{j-k}_Q
\big \|_{L^2}^p\Big)^{1/p}
\endalign
$$
and from here on the proof goes exactly as for Proposition 4.1.\qed
\enddemo

\proclaim{Proposition 4.3}For $1\le p\le 2$
$$
\|S_{m,\s} F\|_{L^p} \lc
2^{-m(d-1)/p'} 2^{-\s(d-1)(1/p-1/2)}
\sup_l\Big(\sum_{j\ge 0}  N(E^{j+l}, 2^{l})
2^{-j(d-1)p/p'}\Big)^{1/p}
\|F\|_{X^p}.
\tag 4.11 
$$
\endproclaim

\demo{\bf Proof} 
We have
(using  Lemma 2.2 (iii) for  the sets $J\in \cJ^k_\mu$
and noting $N(J,2^{k-\mu-m}) \approx \card(J)$)
$$
\align
\|S_{m,\s}F\|_{L^p}&\lc
\Big\|\Big(\sum_{k,\mu}\sum\Sb J\in \cJ^k_\mu
\\ b_J-a_J> 2^{k-\mu-m}\endSb
\big[
\sup_{t\in J}|\cA_t^{\mu+m} \sum_{L(Q)=k-\mu-m+\s}
F_Q^{\mu+m-k}
|\big]^p\Big)^{1/p}\Big\|_{L^p}
\\
&\lc \Big(\sum_{k,\mu}\sum_{J\in \cJ^k_{\mu}}
\big[ 2^{-\s(d-1)(1/p-1/2)} 2^{-(\mu+m)(d-1)/p'}\card(J)^{1/p}
\big]^p
\\&\qquad\qquad\times
\Big(\sum_{L(Q)=k-\mu-m+\s}|Q|^{1-p/2}\|
F_Q^{\mu+m-k}\|_{L^2}^p\Big)\Big)^{1/p}.
\endalign
$$
Now by (1.10)
the latter expression  is estimated
by $2^{-\s(d-1)(1/p-1/2)}2^{-m(d-1)/p'}$ times the quantity
$$
\Big(\sum_{k,\mu}N(E^k, 2^{k-\mu}) 2^{-\mu(d-1)p/p'}
\sum_{L(Q)=k-\mu-m+\s}|Q|^{1-p/2}\|
F_Q^{\mu+m-k}\|_{L^2}^p\Big)^{1/p}
$$
which is bounded by
$$\sup_\ell\Big\{\sum_\mu N(E^{\ell+\mu}, 2^\ell)2^{-\mu(d-1)p/p'}\Big\}^{1/p}
\Big(\sum_l  \sum_{L(Q)=l-m+\s}|Q|^{1-p/2}\|
F_Q^{m-l}\|_{L^2}^p\Big)^{1/p}.
$$
This gives the claimed estimate.\qed
\enddemo

\proclaim{Proposition 4.4}
Suppose that
$$\sup_k
\sum_{j\ge 0}   2^{-j(d-1)p/p'}
N(E^{j+k},2^{k})\le C_1^p.
$$
Then for  $1\le p\le 2$ we have the inequality
$$
\|M_{\ell,\sigma} F\|_{L^p}\,\lc\, C_1
2^{-\sigma (d-1)(1/p-1/2)}
 2^{-\ell(1-d/p')}\|F\|_{X^p}.
\tag 4.12
$$
\endproclaim

\demo{\bf Proof}
This is a small (but crucial) variation of the proof of Proposition 4.3.
We have by Lemma 2.3 part (ii),
$$
\align
\|M_{\ell,\sigma } F\|_{L^p}
&\le \Big(
\sum_k
\sum_{\mu>\ell + \s} \sum_{J\in\cJ^{k,\ell}_\mu}\big\|\sup_{t\in J}
| \cA_t^{\mu-\ell}G_\sigma ^{\mu-\ell-k}(F)
|\big\|_{L^p}^p\Big)^{1/p}
\\
&\lc 2^{-\sigma(d-1)(1/p-1/2)}\, \times
\\&\quad\Big(
\sum_k
\sum_{\mu>\ell} \sum_{J\in\cJ^{k,\ell}_\mu}
\card(J)
2^{-(\mu-\ell)(d-1)p/p'}  2^{-\ell}
\sum_{Q:L(Q)=\sigma-\mu+\ell+k} |Q|^{1-p/2}
\|F_Q^{\mu-\ell-k}\|_{L^2}^p\Big)^{1/p}
\endalign
$$
and this expression by (1.10) is controlled by
$ 2^{-\ell(1-d/p')} 2^{-\sigma(d-1)(1/p-1/2)}$ times the expression 
$$
\Big(
\sum_{n}
\sum_{Q:L(Q)=\sigma+n+\ell} |Q|^{1-p/2}\|F_Q^{-n-\ell}\|_2^p
\sum_{\mu>\ell}
 N(E^{n+\mu}, 2^{n})2^{-\mu(d-1)p/p'}
\Big)^{1/p}
$$
which is
$$\lc
\sup_n
\Big(
\sum_{\mu>\ell}
 N(E^{n+\mu}, 2^{n})2^{-\mu(d-1)p/p'}\Big)^{1/p}
\|F\|_{X^p } .
$$
Thus (4.12) follows.
\qed
\enddemo

\demo{\bf Proof of Theorem I} Immediate by Propositions 3.3, 4.2, 4.3
and 4.4.\enddemo

\head{\bf 5. Weak type ($\boldkey p,\boldkey p$) estimates}\endhead
In this section we shall  mostly  assume that $p<d/(d-1)$ and 
$$
\sup_{k}\sup_{j\ge 0}
N(E^k, 2^{k-j})^{1/p}2^{-j(d-1)/p'}\le C_0.
\tag 5.1
$$
Some statements however will extend to the limiting case $p=d/(d-1)$.

The proof of  Theorem  II follows from 
Proposition 3.3, (4.6), (4.9),
Proposition 4.2 and 
estimates for the operators $S_{m,\s}$ and $M_{\ell,\s}$, stated in
the
following Propositions 5.1 and 5.2. 

\proclaim{Proposition 5.1}
Let $S_{m,\s} F$ be as in (4.4).
Suppose that  $1< p\le d/(d-1)$ if $d=3$, and $1<p<d/(d-1)=2$ if $d=2$, and assume that
(5.1) is valid.
Then there is $\eps=\eps(p)>0$ so that for all $\s, m \ge 0$
 $$
\|S_{m,\sigma } F\|_{L^{p,\infty}}\lc 2^{-\eps(\sigma +m)}\|F\|_{X^p}.
\tag 5.2
$$\endproclaim


\proclaim{Proposition 5.2}
Let $M_{\ell,\sigma } F$ be as in (4.8).
Suppose that  $1< p<d/(d-1)$ and that (5.1) holds. Then  there is
$\eps=\eps(p)>0$ so that
for  $\sigma , \ell \ge 0$
$$
\|M_{\ell,\sigma } F\|_{L^{p,\infty}}\lc 2^{-\eps(\sigma+\ell) }
\|F\|_{X^p} .
\tag 5.3
$$\endproclaim

\demo{\bf Proof of Proposition 5.1}
We have to show that for every $\alpha>0$
$$
\meas(\{x:|S_{m,\sigma } F(x)|\ge 3 \alpha\})\lc 2^{-\eps(m+\sigma )p}
\alpha^{-p} \|F\|_{X^p}^p.
\tag 5.4
$$
Now fix $\alpha>0$ and let 
$$c_Q=|Q|^{1/p-1/2}\Big(\sum_{l:L(Q)+l\ge 0}\|F^l_Q\|_{L^2}^2\Big)^{1/2};
\tag 5.5$$
so that $\sum c_Q^p = \|F\|_{X^p}^p$.
Fix a small $\eps_0>0$ to be chosen later.  
We  divide up the dyadic cubes into two 
families;
$$
\cG = \{ Q :  c_Q^p\frac{1}{|Q|}\le
2^{\eps_0 (\sigma +m) p} \alpha^p \},
\tag 5.6
$$
and  complementary family $\Gamma$, so that $\{Q\}= \cG \cup \Gamma$ and
 $\cG\cap \Gamma=\emptyset$.
Define
$$\aligned
\cG(F)&= \{F^l_Q\}\Sb L(Q)+l\ge 0\\ Q\in \cG\endSb
\\
\cB(F)&= \{F^l_Q\}\Sb L(Q)+l\ge 0\\ Q\in \Gamma\endSb.
\endaligned
$$

For $S_{m,\s}\cG(F)$ we use a straightforward $L^{2}$ estimate.
>From 
Lemma 2.2 (iii) (with $\cE = J \in \cJ^k_{\mu}$), (1.10) and
(5.1) we deduce
$$
\align
\|S_{m,\sigma } \cG(F)\|_{L^{2}}^{2} &\le
\Big\|\Big(\sum_{k} \sum\Sb \mu> 0:\\ \mu + m\ge \s\endSb
\sum\Sb J\in \cJ^k_\mu
\\ b_J-a_J> 2^{k-\mu-m}\endSb
 \sup_{t\in J}\big|\cA_t^{\mu+m}
[\sum \Sb Q\in \cG \\ L(Q)= k-\mu-m+\sigma\endSb 
F^{\mu+m-k}_Q]\big|^{2}\Big)^{1/{2}}\Big\|_{L^{2}}^{2}
\\
&\lc \sum_{k}\sum\Sb \mu> 0:\\ \mu\ge \s - m\endSb 
N(E^k,2^{k-\mu}) 2^{-(\mu+m)(d-1)}
\sum \Sb Q\in \cG\\ L(Q)= k-\mu-m+\sigma \endSb 
\| F^{\mu+m-k}_Q\|_{L^{2}}^{2}
\\
&\lc 2^{-m(d-1)}\sum_{k}\sum \Sb \mu> 0:\\ \mu\ge \s - m\endSb 
2^{-\mu (d-1)(2-p)}
\sum \Sb Q\in \cG\\ L(Q)= k-\mu-m+\sigma \endSb
\| F^{\mu+m-k}_Q\|_{L^{2}}^{2} .
\tag 5.7
\endalign
$$
>From (5.5) and (5.6) we have for $Q\in \cG$ 
$$\| F^{\mu + m-k}_Q \|_{L^2} \le c_Q / |Q|^{1/p - 1/2}
\le 2^{\eps_0 (\s + m)} |Q|^{1/2} \alpha .
\tag 5.8
$$
By \v Ceby\v sev's inequality and (5.7), (5.8) we obtain
$$\align
\meas(\{ & x:\, |S_{m,\sigma }\cG(F)(x)|>\alpha\})\\
&\le
\alpha^{-2}\|S_{m,\sigma }\cG(F)\|_{L^{2}}^{2}
\\
&\lc\alpha^{-2} 2^{-m(d-1)} \sum_{k}\sum\Sb \mu> 0:\\ \mu\ge \s - m\endSb 
 2^{-\mu(d-1)(2 - p)} \sum\Sb Q\in \cG\\L(Q)=k-\mu-m+\sigma \endSb
\|F_Q^{\mu+m-k}\|_{L^2}^{p} \|F_Q^{\mu+m-k}\|_{L^2}^{2 - p}
\\
&\lc
\alpha^{-p} 2^{\eps_0 \s (2 - p)} 
2^{-m [(d-1)-\eps_0 (2-p)]}
\sum_k \sum \Sb \mu> 0:\\ \mu\ge \s - m\endSb 2^{-\mu(d-1)(2-p)}
\sum\Sb Q\in \cG\endSb |Q|^{1-p/2} \big\|F_Q^{\sigma-L(Q)}\big\|_{L^2}^{p}
\\
&\lc
\alpha^{-p} 2^{-\eps \s} 
2^{-\eps m} \|F\|_{X^p}^p 
\tag 5.9
\endalign
$$
for some $\eps > 0$ if $\eps_0 >0$ is small enough. 

We now concentrate on the family $\Gamma$ of dyadic cubes which do not
belong to $\cG$.
Define
$$A(Q,\tau)\equiv
A_{\alpha,\sigma ,m}(Q,\tau):=
2^{(\sigma +m)\eps_0 p} 
\alpha^p 2^{\tau (d-1)p} 2^{L(Q)[\frac 1p-\frac{d-1}{p'}]p} ;
\tag 5.10
$$
note that $\tau\mapsto A(Q,\tau)$ defines an increasing
unbounded sequence for $\tau\ge L(Q)$ and in particular 
$$A(Q,L(Q))=
2^{(\sigma +m)\eps_0 p} \alpha^p |Q|
\tag 5.11
$$
so that for every $Q\in \Gamma$, $c_Q^p > A(Q,L(Q))$.

\definition{Definition}
For every
$Q \in \Gamma$ we define $\tau(Q)$ to be the smallest integer
$\tau>L(Q)$ so that
$A(Q,\tau)\ge c_Q^p$.

\enddefinition

For each $Q\in \Gamma$ we then define $k(Q,\gamma) =  
(L(Q)+\gamma \tau(Q))/(\gamma +1)$ and
$$W(Q)= \bigcup_{k(Q,\gamma)<k\le \tau(Q)}
\bigcup_{I\in\fI_{L(Q)}(E^k)}
\{x\in \Bbb R^d: \big||x-x_Q|-r_I\big|\le 2^{L(Q)+4} 2^{(\tau(Q)-k)\gamma}\}
\tag 5.12
$$
where $\gamma <(d-1)p$ and note that
$$
\align
\meas(W(Q))&\lc \sum_{k\le \tau(Q)}
N(E^k, 2^{k-(k-L(Q))}) 2^{L(Q)+k(d-1)} 2^{(\tau(Q)-k)\gamma }
\\
&\lc \sum_{k\le \tau(Q)}
2^{(\tau(Q)-k)\gamma } 2^{(k-L(Q))\frac{d-1}{p'}p} 2^{L(Q)+k(d-1)}
\\
&\lc  2^{\tau(Q) (d-1)p} 2^{L(Q)[\frac 1p-\frac{d-1}{p'}]p} .
\tag 5.13
\endalign
$$

Let
$$
\cW=\bigcup\Sb Q\in \Gamma\endSb
\big(\{x\in {\Bbb R}^d : |x-x_Q| \le 2^{k(Q,\gamma) + 4}
\} \cup W(Q)\big).
$$
By (5.10), (5.13) and the definition of $\tau(Q)$ 
$$
\align
\meas(\cW) \lc 
\sum\Sb Q\in \Gamma\endSb [2^{k(Q,\gamma) d}
&+ \meas(
W(Q))] \lc 2^{\tau(Q) (d-1)p} 2^{L(Q)[\frac 1p-\frac{d-1}{p'}]p}
\\
&\lc
\sum\Sb Q\in \Gamma\endSb
 2^{-(\sigma +m)\eps_0 p} \alpha^{-p}
A(Q,\tau(Q))
\\
&\lc  2^{-(\sigma +m)\eps_0 p} \alpha^{-p}  \sum\Sb Q\in \Gamma
\endSb
c_Q^p
\lc
2^{-(\sigma +m)\eps_0 p} \alpha^{-p} \|F\|_{X^p}^p.
\endalign
$$

It remains to be shown that
$$
\meas(\{x\notin \cW: S_{m,\sigma}(\cB(F))>2\alpha\})
\lc 2^{-(\sigma +m)\eps_0 p}
 \alpha^{-p} \|F\|_{X^p}^p.
\tag 5.14
$$
We  split
$S_{m,\sigma}(\cB(F))=\sum_{s=-\infty}^\infty {I_s}$ where
$$
I_s=\sup_{k}\sup_{\mu + m - \s \ge \max\{s,0\}}
\sup\Sb
J\in\cJ^k_\mu\\ b_J-a_J\ge 2^{k-\mu-m}\endSb\sup_{t\in  J}
\big|\cA_t^{\mu+m}\big[
\sum\Sb Q\in \Ga\\
L(Q)= k-\mu-m+\sigma \\ \tau(Q)= k-s\endSb
 F^{\mu+m-k}_Q\big]\big|.
$$

We shall prove  
$$
\|I_s\|_{L^{2}}^{2} \lc
2^{-s(d-1)(2 -p)}
2^{-\s (2-p)(d-1-\eps_0) } 
2^{- m [(d-1)(p-1) - \eps_0 (2 - p)]}\alpha^{2-p}
\|F\|_{X^p}^p,
\quad s\ge 0,
\tag 5.15
$$
and 
$$\|I_s\|_{L^p(\Bbb R^d\setminus\cW)}^p\le C_M
2^{-M(\sigma+\gamma|s|)(2-p)}   2^{-m(d-1)(p-1)}
\| F\|_{X^p}^p, \qquad  s\le 0.
\tag 5.16
$$
Note that for $\eps_0 > 0$ small enough
inequalities (5.15) and (5.16) imply (5.14) since
$$
\align
\meas(\{x\notin \cW:\, & S_{m,\sigma}(\cB(F))>2\alpha\})
\\
&\le
\alpha^{-2}\Big\|\sum_{s\ge 0}I_s\Big\|_{L^{2}}^{2}
+
\alpha^{-p}\Big\|\sum_{s<0}I_s\Big\|_{L^p(\Bbb R^d\setminus \cW)}^p
\\&\lc 2^{-\eps(\sigma+m)p}\alpha^{-p}\|F\|^p_{X^p}
\tag 5.17
\endalign
$$
for suitable $\eps=\eps(p)>0$.

\demo{Proof of (5.15)}
We use Lemma 2.2 (iii) for $\cE = J \in \cJ^k_{\mu}$, (1.10) and (5.1)
to obtain 
$$\align
\|I_s\|_{L^{2}}^{2}&\lc
\sum_{k,\mu,J} 
\big\|\sup_{t\in J} \cA_t^{\mu+m}\big[
\sum\Sb Q\in \Ga\\
 \tau(Q)= k-s\\
L(Q)=
k-\mu-m+\sigma
\endSb
 F^{\mu+m-k}_Q\big]\big\|_{L^{2}}^{2}
\\&\lc
\sum_{k,\mu,J} \card(J) 2^{-(\mu+m)(d-1)}
\sum\Sb Q\in \Ga\\ \tau(Q)=k-s\\ L(Q)=k-\mu-m+\sigma\endSb
\|F_Q^{\mu+m-k}\|_{L^2}^{2}
\\&\lc
2^{-m(d-1)}
\sum_{k,\mu} 2^{-\mu(d-1)(2 -p)}
\sum\Sb Q\in \Ga\\ \tau(Q)=k-s\\ L(Q)=k-\mu-m+\sigma\endSb
\|F_Q^{\mu+m-k}\|_{L^2}^{2} .
\endalign
$$
As $k=\tau(Q)+s$ and $\mu=\tau(Q)-L(Q)+s+\s-m$
this  inequality can be rewritten as
$$
\|I_s\|_{L^2}^2\lc 2^{-m(d-1)}\sum_{Q\in \Ga}
2^{(\tau(Q)-L(Q)+s+\s-m)(d-1)(p-2)}
\|F_Q^{\s-L(Q)}\|_{L^2}^2.
\tag 5.18
$$
Now we use that for $Q\in \Ga$
$$\align
&\|F_Q^{\s-L(Q)}\|_{L^2}^{2-p}
\le (c_Q|Q|^{1/p-1/2})^{2-p} \le
\big(2^{-d(1/p-1/2)L(Q)} A(Q,\tau(Q))\big)^{2-p}
\\
&\lc \Big[\alpha 2^{\eps_0(\s+m)} 2^{(d-1)\tau(Q)} 
2^{L(Q)(\frac 1p -\frac {d-1}{p'} -d(\frac 1p-\frac 12))}\Big]^{2-p}
\endalign
$$ and combine this with (5.18) to obtain after a little algebra
$$
\|I_s\|_{L^2}^2
\lc 2^{m((d-1)(1-p)+\eps_0(2-p))} 2^{-\s(d-1-\eps_0)(2-p)}
 2^{-s(d-1)(2-p)} \alpha^{2-p} \sum_Q|Q|^{1-p/2}\|F_Q\|_{L^2}^p
$$
which is the desired bound.
\enddemo

\demo{Proof of (5.16)}
We use the 
estimate away from the exceptional set in 
Lemma 2.2 (iv),  with $\eta= 2^{|s|\gamma}$ ($\gamma<(d-1)p$) 
and $s=k-\tau(Q)$.
Then
$$\align
\|I_s\|_{L^p(\Bbb R^d\setminus\cW)}^p
&\lc
\sum_{k,\mu} 
\big\|\sup_{t\in J} \cA_t^{\mu+m}\big[
\sum\Sb
\tau(Q)= k-s\\L(Q)=  k-\mu-m+\sigma
\endSb
 F^{\mu+m-k}_Q\big]\big\|_{L^p(\Bbb R^d\setminus \cW)}^p
\\&\lc
\sum_{k,\mu} N(E^k,2^{k-\mu}) 2^{-(\mu+m)(d-1)p/p'}
2^{-(\sigma+\gamma |s|)M(2-p)}
\sum\Sb\tau(Q)=k-s\\ L(Q)=k-\mu-m+\sigma\endSb
|Q|^{1-p/2}\|F_Q^{\mu+m-k}\|_{L^2}^p
\\&\lc
 2^{-m(d-1)p/p'}
2^{-\sigma M(2-p)} 2^{-|s|\gamma (2-p)}
\|F\|_{X^p}^p.\qed
\endalign
$$
\enddemo
\enddemo

\subheading{Proof of Proposition 5.2}

This is similar to the proof of Proposition 5.1.
We have to show that for every $\alpha>0$
$$
\meas(\{x:|M_{\ell,\sigma } F(x)|\ge 3 \alpha\})\lc 2^{-\eps(\ell+\sigma )p}
\alpha^{-p} \|F\|_{X^p}^p.
\tag 5.19
$$

We indicate the changes in the proof of Proposition 5.1. Of course we
systematically replace  $S_{m,\sigma}$ by $M_{\ell,\sigma}$.
The definition (5.6) is  the same except that $2^{\eps_0 mp}$ has
to be replaced by $2^{\eps_0 \ell p}$; then the arguments up to
(5.9) are similar; we have to use Lemma 2.3 
(ii) instead of Lemma 2.2 (iii).
Similarly the definition (5.10) is changed to
$$A(Q,\tau)\equiv A_{\alpha,\sigma ,\ell }(Q,\tau):=
2^{(\sigma +\ell)\eps_0 p} 
\alpha^p 2^{\tau (d-1)p} 2^{L(Q)[\frac 1p-\frac{d-1}{p'}]p};
$$
and the further arguments up to (5.14) have obvious analogues.
In the definition of $A(Q,\tau)$ we shall need to take $\eps_0$ so that
$\eps_0(2-p)<1-(d-1)(p-1)$ which is possible since $p<d/(d-1)$.

We then split $M_{\ell,\sigma}(\cB(F))=\sum II_s$
where
$$\aligned II_s\,&=\,
\sup_k
\sup_{\mu\ge \ell+ \s}\sup_{J\in\cJ^{k,\ell}_\mu}\sup_{t\in J}
\big|\cA_t^{\mu-\ell}G_{\sigma,s}^{\mu-\ell-k}(F)
\big|
\\&\text{ and }
G_{\sigma,s}^{\mu-\ell-k} F :=
\sum\Sb Q: L(Q)=k-\mu+\ell+\sigma\\ \tau(Q)=k-s\endSb F^{\mu-\ell-k}_Q.
\endaligned
\tag 5.20
$$

The inequalities (5.15) and (5.16) are replaced by

$$
\|II_s\|_{L^{2}}^{2} \lc
2^{-s(d-1)(2 -p)}
2^{-\s (2-p)(d-1-\eps_0)}
2^{- \ell [(1-(d-1)(p-1)) - \eps_0 (2 - p)]}\alpha^{2-p}
\|F\|_{X^p}^p,
\quad s\ge 0,
\tag 5.21
$$
and
$$
\|II_s\|_{L^p(\bbR^d\setminus\cW)}^p \lc 2^{- \ell(1-\frac{d}{p'})p}
2^{-s(d-1)(1-p/2)}
2^{-(\sigma+|s|\gamma)Mp}
\alpha^{2-p}
\| F\|_{X^p}^p, \qquad s\le 0,
\tag 5.22
$$
from which we can  as before  conclude the assertion of the proposition.

%
%

\demo{Proof of (5.21) and (5.22)} We prove (5.21) and use 
Lemma 2.3 to estimate
$$
\align
\|II_s\|_{L^2}^2
&\le \sum_k\sum_{\mu\ge \ell+\s}
\sum_{J\in \cJ^{k,\ell}_\mu}\big\|\sup_{t\in J}|\cA_t^{\mu-\ell} 
G^{\mu-\ell-k}_{\s,s}(F)\big\|_{L^2}^2
\\
&\le \sum_k\sum_{\mu\ge \ell+\s}
\sum_{J\in \cJ^{k,\ell}_\mu} \card(J) 2^{-(\mu-\ell)(d-2)}2^{-\mu}
\sum \Sb Q\in \Ga\\ 
L(Q)=k-\mu+\ell+\s\\ \tau(Q)=k-s
\endSb
\|F_Q^{\mu-\ell-k}\|_{L^2}^2
\endalign
$$
Now 
$\sum_{J\in \cJ^{k,\ell}_\mu}\card J\lc N(E^k, 2^{k-\mu})\lc 
2^{\mu(d-1)(p-1)}$ by assumption (1.10) and ($\cC_{p,\infty}$). 
We also observe that 
$\mu=\tau(Q)-L(Q)+s+\s+\ell$ in the above sum and thus we obtain
$$\align
\|II_s\|_{L^2}^2&\lc 2^{\ell(d-2)}
\sum_k\sum_{\mu\ge \ell +\s} 2^{\mu(d-1)(p-2)}
\sum \Sb Q\in \Ga\\ 
L(Q)=k-\mu+\ell+\s\\ \tau(Q)=k-s
\endSb
\|F_Q^{\mu-\ell-k}\|_{L^2}^2
\\
&\lc 2^{\ell(d-2)}
\sum\Sb Q\in \Ga\\ \tau(Q)- L(Q)\ge -s\endSb 
\|F_Q^{\s-L(Q)}\|_{L^2}^2
2^{(\tau(Q)-L(Q)+s+\s+\ell)(d-1)(p-2)}.
\endalign
$$
Now as before
$\|F_Q^{\s-L(Q)}\|_{L^2}^{2-p}
\lc \big[\alpha 2^{\eps_0(\s+\ell)} 2^{(d-1)\tau(Q)} 
2^{L(Q)(\frac 1p -\frac {d-1}{p'} -d(\frac 1p-\frac 12))}\big]^{2-p}
$
and after doing the algebra we arrive at
$$
\|II_s\|_2^2\lc  \alpha^{2-p}
2^{\ell((d-1)p-d+\eps_0(2-p))} 2^{-\s(d-1-\eps_0)(2-p)} 2^{-s(d-1)(2-p)}
\sum_{Q\in \Ga} |Q|^{1-p/2}\|F_Q^{\s-L(Q)}\|_{L^2}^p
$$
which is what we were aiming for.

Similarly, the proof of (5.22) is analogous to the proof of (5.16).
\qed\enddemo


\head{\bf 6. $\boldkey L^{\boldkey p}$ estimates in the limiting case}\endhead

We assume throughout this section that that the
regularity condition $(\widetilde \cR_{p_d})$, $p_d=d/(d-1)$, is satisfied.
We first give a reformulation of the Carleson-measure condition.

\proclaim{Lemma 6.1} Suppose that the Carleson measure condition $(\widetilde 
\cC_{p_d})$ holds.
Then the measure 
$$\sum_{k\in \bbZ}\sum_{\mu\ge 0} \delta_{k,\mu} 
\sum_{J\in \cJ^k_\mu} \card(J) 2^{-\mu}(1+\mu)^{\frac d{d-1}}$$
is also a Carleson measure.
\endproclaim

\demo{\bf Proof} 
We first observe that
$$
N(E^k, 2^{k-j}) 2^{k-j}\approx \big|\{t\in[2^k, 2^{k+1}):\dist(t, E^k)\le 2^{k-j}\}\big|
$$
and thus 
$$N(E^k, 2^{k-j}) 2^{-j} \le C N(E^k, 2^{k-j'}) 2^{-j'} \text{ if } j'\le j.
\tag 6.1
$$

Let $I$ be an interval of length $>1$ and $I^*$ the interval with same midpoint and double
length.
Then
$$\align
&\sum_{(k,\mu)\in T(I)}\sum_{J\in \cJ^k_\mu}\card( J )2^{-\mu}(1+\mu)^{d/(d-1)}
\\
&\lc
\sum_{k\in I}\sum_{s=0}^{1+\log_2|I|}
2^{sd/(d-1)}\sum_{2^{s-1}\le\mu<2^{s}} 2^{-\mu} \sum_{J\in \cJ^k_{\mu}}\card( J) 
\\
&\lc
\sum_{k\in I}\sum_{s=0}^{1+\log_2|I|}
2^{sd/(d-1)}
N(E^k, 2^{k-2^{s-1}}) 2^{-2^{s-1}}
\\
&\lc
\sum_{k\in I}\sum_{s=0}^{1+\log_2|I|}
\sum_{2^{s-2}\le\mu<2^{s-1}}
N(E^k, 2^{k-\mu}) 2^{-\mu}(1+\mu)^{1/(d-1)}
\\
&\lc
\sum_{(k,\mu)\in T(I^*)}
N(E^k, 2^{k-\mu}) 2^{-\mu}(1+\mu)^{1/(d-1)}.
\endalign
$$

Here we have used the regularity assumption (1.11) for the second 
inequality  and (6.1) for the third
inequality.\qed

\enddemo

The following is an even more elementary observation.

\proclaim{Lemma 6.2} Suppose that the Carleson measure condition $(\widetilde \cC_{p_d})$ holds.
Then
$$\sup_\nu \sum_\mu \sum_{J\in\cJ^{\mu+\nu}_\mu}\card( J) 2^{-\mu} \le C.
\tag 6.2
$$
\endproclaim

\demo{\bf Proof}
Let $I_s(r)=\{x:|x-r|\le 2^s\}$. Then
$$
\align
\sum_\mu \sum_{J\in\cJ^{\mu+r}_\mu}\card( J) 2^{-\mu}
&\lc
\sum_{s=0}^\infty 2^{-sd/(d-1)}\sum_{0\le \mu\le 2^s}
\sum_{J\in\cJ^{\mu+r}_\mu}\card (J) 2^{-\mu}(1+\mu)^{d/(d-1)}
\\
&\lc\sum_{s=0}^\infty 2^{-s/(d-1)} \frac{1}{|I_s(r)|}\sum_{(k,\mu)\in T(I_s(r))}
\sum_{J\in\cJ^{k}_\mu }\card( J) 2^{-\mu} (1+\mu)^{d/(d-1)}
\endalign
$$
and the last expression is bounded by Lemma  6.1.\qed
\enddemo

The following Carleson-measure estimate is a standard consequence of the $L^p$ boundedness of the Hardy-Littlewood maximal operator, for the proof see 
\cite{14, ch. II.2}.
\proclaim{Lemma 6.3}
Suppose the doubly indexed nonnegative  sequence $\{\omega_{k,\mu}, (k,\mu)\in \bbZ\times\bbZ^+\}$
satisfies
$$
\sup_{|I|\ge 1} \frac{1}{|I|}\sum_{(k,\mu)\in T(I)}\omega_{k,\mu}\le A^p
$$
i.e. $ \sum\omega_{k,\mu}\delta_{k,\mu}$ is a Carleson measure. Then for $\{a_k\}\in\ell^p$, $p>1$
$$
\Big(\sum_{k,\mu}\omega_{k,\mu}\Big[\frac1{1+\mu}\sum_{|j|\le \mu} |a_{k+j}|\Big]^p\Big)^{1/p}
\le C_p A \Big(\sum_k |a_k|^p\Big)^{1/p}.
$$
\endproclaim

We now turn to the
\subheading{$\boldkey L^{\boldkey 2}$ estimates in two dimensions}
We are concerned with the $L^2(\Bbb R^2)$ estimates in  Theorem III.
The claim is a consequence of the following  estimates:

$$
\Big\|\sup_{k,\mu}\sup_{J\in\cJ^k_{\mu}}
\sup_{t\in J}\big|\sum\Sb  j\le \mu\\ b_J-a_J> 2^{k-j}\endSb
\cA_t^j f\big|\Big\|_{L^2}\lc 
\|f\|_{L^2}
\tag 6.3
$$
and, for $m\ge 0$,
$$
\Big\|\sup_{k,\mu}\sup\Sb J\in\cJ^k_{\mu}\\ b_J-a_j>2^{k-\mu-m}\endSb
\sup_{t\in J}\big| \cA_t^{\mu+m} f\big|\Big\|_{L^2}\lc  2^{-m/2}
\|f\|_{L^2}
\tag 6.4
$$
and finally
$$
\Big\|\sup_{k}\sup_{J\in\cJ^{k}}
\sup_{t\in J}\big|\sum\Sb j:\\ b_J-a_J\le 2^{k-j}\endSb
\cA_t^j f\big|\Big\|_{L^2}\lc 
\|f\|_{L^2}.
\tag 6.5
$$

To prove (6.3) we use  Lemma 2.3 to see that the left side  is dominated by
$$
\align
&\Big(\sum_{k,\mu}\sum_{J\in\cJ^k_{\mu}}\big[\sum\Sb j\le \mu\\ b_J-a_J>2^{k-j}\endSb
 \|\sup_{t\in J}
|\cA_t^j f| \|_{L^2}\big]^2\Big)^{1/2}
\\
&\lc\Big(\sum_{k,\mu}\sum_{J\in\cJ^k_{\mu}}
\card( J) 2^{-\mu}(1+\mu)^2\Big[\frac{1}{1+\mu}\sum_{j\le \mu}
\|P^{j-k} f\|_{L^2}\Big]^2
\Big)^{1/2}
\endalign
$$
and by Lemma 6.3 and 6.1 the last expression is controlled  by 
$$\Big(\sum_{k\in \bbZ}\|P^k f\|_{L^2}^2\Big)^{1/2}\lc \|f\|_{L^2}.$$

Concerning  (6.4) we use Lemma 2.2 and bound the left side by
$$
\align
&\Big(\sum_{k,\mu}\sum\Sb J\in\cJ^k_{\mu}\\b_J-a_J>2^{k-j}\endSb\big\|
\sup_{t\in J}\big| \cA_t^{\mu+m} f\big| \big\|_{L^2}^2\Big)^{1/2}
\\
&\lc\Big(\sum_{k,\mu}\sum_{J\in\cJ^{k}_{\mu}} 
N(J, 2^{k-m-\mu})
2^{-(\mu+m)}\|P^{\mu+m-k}f\|_{L^2}^2\Big)^{1/2}
\\
&\lc 2^{-m/2}\sup_{l\in \bbZ} \Big(\sum_\mu
\sum_{J\in\cJ^{\mu+m-l}_{\mu}} 
\card (J)
2^{-\mu}\Big)^{1/2}
\Big(\sum_k\|P^{k}f\|_{L^2}^2\Big)^{1/2}
\endalign
$$
and by Lemma 6.2 the last expression is $\lc 2^{-m/2} \|f\|_{L^2}$.

Finally (6.5) holds in view of the assumption (1.9);
\cf. the argument in the proof of Proposition 4.1. We shall not repeat the details.
\qed

\subheading {$\boldkey X^{\boldkey p}$ estimates and the proof of Theorem III }
We use a similar decomposition as in \S4 however instead of considering the 
maximal operators $M_{\ell,\sigma}$ we shall not decompose in $\ell$ and work with $M_\s$ in (4.5) directly.
We shall prove
$$\|M_\s F\|_{L^{p_d}}\lc 2^{-\s(d-1)(1/p_d-1/2)}\|F\|_{X^{p_d}} .
\tag 6.6
$$
This together with already proved estimates in \S4 implies the statement of Theorem III.

We argue as before and set 
$a_l=(\sum_{L(Q)=\s-l}|Q|^{1-{p_d}/2}\|F_Q^l\|_{L^2}^{p_d})^{1/{p_d}}$.
Using Lemma 2.3 (ii) we get
$$\align
\|M_\s F\|_{L^{p_d}}
&\lc
\Big(
\sum_k\sum_\mu\sum\Sb J\in\cJ^k_\mu\endSb\Big[
\sum\Sb j\le \mu\\ b_J-a_J>2^{k-j}\endSb
\big\|\sup_{t\in J}|
\cA_t^j G^{j-k}_\s(F)|\big\|_{L^{p_d}}
\Big]^{p_d}\Big)^{1/{p_d}}
\\
&\lc
\Big(
\sum_k\sum_\mu\sum_{J\in\cJ^k_\mu}
\Big[
\sum_{10< j\le\mu} 
 2^{-\s(d-1)(1/{p_d}-1/2)}
\card(J)^{1/{p_d}}2^{-\mu/{p_d}}
a_{j-k}\Big]^{p_d}\Big)^{1/{p_d}}
\\
&\lc
 2^{-\s(d-1)(1/{p_d}-1/2)}
\Big(
\sum_k\sum_\mu\sum_{J\in\cJ^k_\mu} \card(J) 2^{-\mu}\mu^{p_d}
\Big(\frac{1}{\mu+1}\sum_{0\le j\le \mu} a_{j-k}\Big)^{{p_d}}\Big)^{1/{p_d}} .
\endalign
$$
By condition $(\widetilde\cC_{p_d})$ and Lemma 6.3 and Lemma 6.1 we obtain 
(6.6).\qed

\head{\bf 7. Weak type ($\boldkey p$, $\boldkey p$) estimates in the limiting case}\endhead

Throughout this section we shall assume that  $d\ge 3$ and
that  the regularity assumption and condition (1.14) hold; thus
$$
\sup_k  2^{-n} N(E^k , 2^{k-n}) \le
C \big( n^{\frac 1{d-1}} \log n\big)^{-1}
\tag 7.1
$$
uniformly in $n\ge 10$. We follow the proof of Theorem II
in \S 5, using the same decompositions except we do not 
decompose $M_{\s}$ in (4.5) further as in the proof of Theorem III.
We recall that  Proposition 5.1 remains valid for the limiting case  $p=p_d$ if $d\ge 3$,
under the weaker condition $(\cC_{p_d,\infty})$. Therefore
the claim in Theorem IV will be a consequence of 

\proclaim{Proposition 7.1} Let $M_{\s}F$ be as in (4.5). Suppose
(7.1) holds. Then there is an $\eps >0$ so that for all
$\s, \alpha \ge 0$,
$$
\meas(\{ x: M_{\s} F(x) > 3 \alpha \}) \lc 2^{-\eps \s p_d} \alpha^{
- p_d} \|F\|_{X^{p_d}}^{p_d} .
\tag 7.2
$$
\endproclaim

\demo{\bf Proof}
As in \S 5 we fix  $\eps_0 >0$ and define $\cG$, $\Gamma$, $A(Q,\tau)$, $\cG (F)$,
$\cB (F)$ and $\cW$ as before except we replace $2^{\eps_0 (\s + m)p}$
with $2^{\eps_0 \s p_d}$. In particular we  have now  for $\tau\ge L(Q)$
$$A(Q,\tau)^{1/p_d}= 2^{\s\eps_0} \alpha 2^{\tau(d-1)}.
$$
We shall have to take $\eps_0$ so that $0<\eps_0(2-p_d)<d-2$.

For $M_{\s} \cG (F)$ we use an $L^2$ estimate. From Lemma 2.3 (ii)
and the regularity assumption
(1.11), we deduce
$$
\align
\| M_{\s} \cG (F)\|_{L^2}^2 &\le \sum_k \sum_{\mu\ge\s}
\sum_{J\in \cJ^k_{\mu}} \Big( \sum_{\s\le j\le\mu} \|\sup_{t\in J}
A^j_t (\sum\Sb Q\in \cG\\L(Q)=k-j+\s\endSb F^{j-k}_Q )\|_{L^2}\Big)^2
\\
&\lc \sum_k \sum_{\mu\ge\s}2^{-\mu} \sum_{J\in \cJ^k_{\mu}} \card(J)
\Big(\sum_{\s\le j\le\mu} 2^{-j(d-2)/2} (
\sum\Sb Q\in \cG\\L(Q)=k-j+\s\endSb \| F^{j-k}_Q\|_{L^2}^2 )^{1/2} \Big)^2
\\
&\lc [\sup_k 2^{-\s}N(E^k , 2^{k-\s})]
\Big(\sum_{\s\le j} 2^{-j(d-2)/2} (
\sum\Sb Q\in \cG\endSb \| F^{\s - L(Q)}_Q\|_{L^2}^2 )^{1/2} \Big)^2
\\
&\lc 2^{\s (d-2)} \sum\Sb Q\in \cG\endSb \| F^{\s - L(Q)}_Q\|_{L^2}^2 .
\tag 7.3
\endalign
$$
For $Q\in\cG$ we have
$$\|F_Q^{\s - L(Q)}\|_{L^2} \le c_Q /|Q|^{1/p_d - 1/2} \le 2^{\eps_0 \s (
2-p_d)}|Q|^{1/2} \alpha$$
and therefore by \v Ceby\v sev's inequality and (7.3),
$$
\align
\meas(\{M_{\s} &(F) > \alpha \}) 
\\
&\le \alpha^{-2} \| M_{\s} (F)\|_{L^2}^2
\\
&\lc 2^{-\s[(d-2)-\eps_0 (2-p_d)]} \alpha^{-p_d}
\sum_Q |Q|^{1-p_d /2} \|F^{\s - L(Q)}_Q \|_{L^2}^{p_d}
\\
&\lc 2^{-\eps \s} \alpha^{-p_d} \|F\|_{X^{p_d}}^{p_d}
\endalign
$$
for some $\eps > 0$.

Furthermore the estimate for the measure of the exceptional set $\cW$
 in \S5  is still valid.
Therefore it remains to be shown that 
$$
\meas(\{ x\notin \cW : M_{\s} \cB (F)(x) > 3\alpha \}) \lc 2^{-\eps \s}
\alpha^{-p_d} \|F\|_{X^{p_d}}^{p_d} 
\tag 7.4
$$
for some  $\eps>0$. We may estimate
$$M_{\s} \cB (F) \le III + \sum_{s=0}^\infty IV_s +\sum_{s=-\infty}^{-1}
 V_s$$ where
$$
\align
III &= \sup_k \sup_{\mu\ge\s}\sup_{J\in\cJ^k_{\mu}} \sup_{t\in J}
\Big| \sum\Sb \s\le j\le \mu\\ b_J - a_J \ge 2^{k-j}\endSb
A^j_t \big[\sum\Sb Q\in \Ga\\
 L(Q)=k-j+\s\\ e^{\eps_1|k-\tau(Q)|}\le \mu\endSb F^{j-k}_Q \big]\Big|
\\
IV_s&=
 \sup_k \sup_{\s\le \mu\le e^{\eps_1s}}
\sup_{J\in\cJ^k_{\mu}} \sup_{t\in J}
\Big| \sum\Sb \s\le j\le \mu\\ b_J - a_J \ge 2^{k-j}\endSb
A^j_t \big[\sum\Sb  Q\in \Ga\\L(Q)=k-j+\s\\ \tau(Q)=k-s\endSb F^{j-k}_Q \big]\Big|, \quad s>0,
\\
V_s &= \sup_k \sup_{\s\le \mu\le e^{\eps_1|s|}}\sup_{J\in\cJ^k_{\mu}} \sup_{t\in J}
\Big| \sum\Sb \s\le j\le \mu\\ b_J - a_J \ge 2^{k-j}\endSb
A^j_t\big[\sum\Sb  Q\in \Ga\\ L(Q)=k-j+\s\\ \tau(Q) = k-s\endSb F^{j-k}_Q \big]\Big|, \quad s<0.
\endalign
$$
Here we may choose $0<\eps_1< (d-2)/2$. We then 
 prove 
$$
\align
&\|III\|_{L^{p_d}}^{p_d} \lc
2^{-\sigma(d-1)(1-p_d/2)} \log (2+\s)
 \|F\|_{X^{p_d}}^{p_d}, \quad  
\tag 7.5
\\
&\|IV_s \|_{L^2}^2 \lc 2^{-\s(d-2-\eps_0(2-p_d))}2^{-s(d-2-2\eps_1)}
\alpha^{2-p_d} \|F\|_{X^{p_d}}^{p_d}, \quad  s\ge 0,
\tag 7.6
\\
&\|V_s\|_{L^{p_d} ({\Bbb R}^d \setminus \cW )}^{p_d} \lc 2^{-M(2-p_d)
(\s + \gamma| s|)} \| F\|_{X^{p_d}}^{p_d}, \qquad s<0.
\tag 7.7
\endalign
$$

(7.4) follows from (7.5), (7.6) and (7.7) in the usual way.
We remark that our assumption (1.14) is needed for (7.5). 
For the error terms (7.6), (7.7) we can get away with just  the 
regularity hypothesis (1.11) and ($\cC_{p_d , \infty}$).

In the proof we shall use arguments that occur in the proof of 
Hardy's inequality
(see \cite{5}).


%

\demo{Proof of (7.5)}
We further split
$III=\sum\Sb 2^n\ge \sigma\endSb III_n$ where
$$
III_n= \sup_k \sup_{\mu\ge 2^n}\sup_{J\in\cJ^k_{\mu}} \sup_{t\in J}
\Big| \sum\Sb 2^{-n}\mu< j\le 2^{-n+1} \mu\\ b_J - a_J \ge 2^{k-j}\endSb
A^j_t \big[\sum\Sb Q\in \Ga\\ L(Q)=k-j+\s\\ e^{\eps_1|k-\tau(Q)|}
\le \mu\endSb F^{j-k}_Q \big]\Big|.
$$
We replace various $\sup$'s by $\ell^{p_d}$ norms and use Lemma (2.3), (ii). We obtain
$$\multline
\|III_n\|_{L^{p_d}}^{p_d}\le \sum_k\sum_{\mu\ge 2^n}\sum_{J\in \cJ^k_\mu}
\Big[
\sum \Sb  2^{-n}\mu< j\le 2^{-n+1} \mu\endSb 
2^{(j-\mu)/p_d}  2^{-j(d-1)/p_d'} \card(J)^{1/p_d}
2^{-\s(d-1)(1/p_d-1/2)}
\\ \times
\Big(\sum\Sb Q\in\Gamma\\ L(Q)=k-j+\s\\|k-\tau(Q)|\le \eps_1^{-1}\log \mu\endSb 
|Q|^{1-p_d/2} \|F_Q^{j-k}\|_{L^2}^{p_d}\Big)^{1/p_d}\Big]^{p_d}.
\endmultline
$$
If we abbreviate 
$$
\aligned
&w_{\mu,k}= 2^{-\mu}\sum_{J\in \cJ^k_\mu}\card J,
\\
&b_{Q,\s}=
|Q|^{1-p_d/2} \|F_Q^{\s-L(Q)}\|_{L^2}^{p_d},
\endaligned
\tag 7.8
$$
this yields
$$\align
&\|III_n\|_{L^{p_d}}^{p_d}
\\
&\lc 2^{-\s(d-1)(1-p_d/2)}
 \sum_k\sum_{\mu\ge 2^n} w_{k,\mu}
\Big[ \sum_{2^{-n}\mu< j\le 2^{-n+1} \mu}
\Big(\sum\Sb Q\in \Ga\\L(Q)=k-j+\s\\|k-\tau(Q)|\le \eps_1^{-1}\log \mu\endSb
 b_{Q,\s}\Big)^{1/p_d}\Big]^{p_d}
\\
&\lc 2^{-\s(d-1)(1-p_d/2)}
\sum_k\sum_{\mu\ge 2^n} w_{k,\mu} \mu^{p_d-1}
\sum_{2^{-n}\mu<j\le 2^{-n+1}\mu} 
\sum\Sb L(Q)=k-j+\s
\\|k-\tau(Q)|\le \eps_1^{-1}\log \mu\endSb b_{Q,\s}
\\
&\lc 2^{-\s(d-1)(1-p_d/2)} 2^{-n(p_d-1)}
\sum_Q
\sum\Sb k: |k-\tau(Q)|\le \eps_1^{-1}\log(2^n(k-L(Q)+\s))\endSb
 \sum\Sb\mu:
2^{n-1}(k-L(Q)+\s)
\le \\ \mu\le 2^{n}(k-L(Q)+\s)
\endSb
 w_{k,\mu} \mu^{p_d-1}
 b_{Q,\s}
\tag 7.9
\endalign
$$
Now by the regularity assumption $(\widetilde R_{p_d})$ and by (1.14) we have
$$
 \align
&\sum\Sb\mu:
2^{n-1}(k-L(Q)+\s)
\le\\ \mu\le 2^{n}(k-L(Q)+\s)
\endSb
 w_{k,\mu} \mu^{p_d-1}
\\
&\qquad\qquad\lc 2^{n(p_d-1)}(k-L(Q)+\s+1)^{p_d-1} N(E^k, 2^{k-2^n(k-L(Q)+\s)}) 2^{-2^n(|k-L(Q)|+\s)}
\\&\qquad\qquad\lc  \big[ \log \log (2^{2^n(|k-L(Q)|+\s)})\big]^{-1}
\endalign
$$
and thus  the expression (7.9) is controlled by
$$\align & 2^{-\s(d-1)(1-p_d/2)} 2^{-n(p_d-1)}
\sum_Q b_{Q,\s}
\sum\Sb k: |k-\tau(Q)|\le \eps_1^{-1}(n+\log(k-L(Q)+\s))\endSb
(1+n+ \log (|k-L(Q)|+\s))^{-1}
\\&
\lc 2^{-\s(d-1)(1-p_d/2)} 2^{-n(p_d-1)}(\log(2+\s)+n)
\sum_Q b_{Q,\s}.
\endalign
$$
Hence
$$
\|III_n\|_{L^{p_d}}^{p_d}\lc 2^{-\s(d-1)(1-p_d/2)}  2^{-n(p_d-1)}(\log(2+\s)+n)
\|F\|_{X_{p_d}^{p_d}}
$$ which yields the asserted bound (7.5).
\enddemo
\demo{Proof of (7.6)}
We estimate $IV_s\le \sum_{2^n>\sigma} IV_{s,n}$ where
$$IV_{s,n}=
 \sup_k \sup_{\s\le \mu\le e^{\eps_1 s}}\sup_{J\in\cJ^k_{\mu}} \sup_{t\in J}
\Big|  \sum \Sb  2^{-n}\mu< j\le 2^{-n+1} \mu
\\ b_J - a_J \ge 2^{k-j}\endSb
\cA^j_t \big[\sum\Sb Q\in \Ga\\L(Q)=k-j+\s\\ k-\tau(Q)=s\endSb 
F^{j-k}_Q \big]\Big| .
$$
We apply  H\"older's inequality for the sum in $j$ and apply Lemma 2.3 to get
$$
\align
\|IV_{s,n}\|^2_{L^2}&\lc
\sum_k\sum_{\s\le \mu\le e^{\eps_1 s}}
\sum_{J\in\cJ^k_\mu} 2^{-n}(\mu+1)
\sum \Sb  2^{-n}\mu< j\le 2^{-n+1} \mu
\\ b_J - a_J \ge 2^{k-j}\endSb\Big\|\sup_{t\in J}
\Big|\cA^j_t \big[\sum\Sb Q\in \Ga\\L(Q)=k-j+\s\\ k-\tau(Q)=s\endSb F^{j-k}_Q \big]
\Big|\Big\|_{L^2}^2
\\
&\lc 2^{-n} e^{\eps_1 s}
\sum_k\sum_{\s\le \mu\le e^{\eps_1 s}}
\sum_{J\in\cJ^k_\mu} \card(J)
\sum \Sb  2^{-n}\mu< j\le 2^{-n+1} \mu\endSb    2^{-j(d-2)} 2^{-\mu}
\sum\Sb Q\in \Ga\\L(Q)=k-j+\s\\ k-\tau(Q)=s\endSb\| F^{j-k}_Q\|_{L^2}^2.
\endalign
$$
Now  we use 
$2^{L(Q) d(1/p_d-1/2)}\|F_Q^{j-k}\|_{L^2}\lc 2^{\eps_0\s} 2^{\tau(Q)(d-1)}\alpha$ and that
$k=\tau(Q)+s$, $j=\tau(Q)-L(Q)+s+\s$  and derive 
$$
\multline
\|IV_{s,n}\|_{L^2}^2\lc 2^{-n}e^{\eps_1 s}
\sum\Sb Q\in \Ga\\ \tau(Q)-L(Q)\ge -s\endSb 
\|F_Q^{\s-L(Q)}\|_{L^2}^{p_d}
2^{\eps_0\s(2-p_d)} 2^{\tau(Q)(d-1)(2-p_d)} \alpha^{2-p_d}
\\
\times
2^{-(\tau(Q)-L(Q)+s+\s)(d-2)}
\Big\{
\sum\Sb \sigma\le \mu\le e^{\eps_1 s}\\
\mu\le 2^n(\tau(Q)-L(Q)+s+\s) \\
\mu\ge 2^{n-1}(\tau(Q)-L(Q)+s+\s) \endSb 2^{-\mu}\sum_{J\in \cJ_\mu^{\tau(Q)+s}}\card J 
\Big\}
\endmultline
$$
The expression $\{ \dots\}$ is $O(1)$ by (1.11). 
We compute that 
$(d-2)-d(1/p_d-1/2)(2-p_d)=d(1-p_d/2)$ and
$(d-1)(2-p_d)=d-2$. Thus the last estimate 
simplifies to 
$$\align
\|IV_{s,n}\|_{L^2}^2&\lc 2^{-n}e^{\eps_1 s}
\sum\Sb Q\in \Ga\\ \tau(Q)-L(Q)\ge -s\endSb  
2^{L(Q)d(1-p_d/2)}\|F_Q^{\s-L(Q)}\|_{L^2}^{p_d}
2^{\eps_0\s(2-p_d)}  \alpha^{2-p_d}
2^{-(s+\s)(d-2)}
\\&\lc \alpha^{2-p_d}
2^{-n} 2^{-\s(d-2-\eps_0(2-p_d))} 2^{-s(d-2-2\eps_1)}
\sum\Sb Q\endSb  |Q|^{1-p_d/2}\|F_Q^{\s-L(Q)}\|_{L^2}^{p_d}
\endalign
$$
which implies (7.6).
\enddemo
\demo{Proof of (7.7)}
This $L^{p_d}$ estimate away from the exceptional set follows  by analogous arguments; 
 Lemma 2.3(iii) is used.
We omit the details.
\enddemo
 This completes
the proof of Proposition 7.1. \qed
\enddemo

\head{\bf 8. Examples and counterexamples}\endhead

We consider a simple  class of sets $E$ to  which
Theorems I-IV can be applied. They satisfy the

\proclaim{8.1. Convexity assumption} For each $k\in \Bbb Z$ the set $E^k$ is given
 by
$\{t_\nu^k\}_{\nu=1}^\infty$ where $t_\nu^k$ is a monotone sequence contained in
$[2^k,2^{k+1}]$, so that the sequence $t_{\nu+1}^k -t_{\nu}^k$ is
also monotone.
\endproclaim

The following lemma shows that if ($\cC_{p,\infty}$) holds for
some $p<d/(d-1)$ and $E$ satisfies the convexity assumption; then it also
satisfies the regularity assumption for all $p>1$.

\proclaim{Lemma 8.1.1} Suppose $E$ satisfies the convexity assumption.
Suppose that for some $\beta>0$ the estimate
$$\sup_k N(E^k, 2^{k-n})\,\le\, C\frac{2^n}{ (1+n)^{\beta}}
\tag 8.1$$
holds uniformly in $k\in\Bbb Z$. Then $E$ satisfies regularity assumption $(\cR_p)$
for all $p>1+[(d-1)(\beta+1)]^{-1}$.
Moreover 
it satisfies regularity assumption ($\widetilde\cR_{d/(d-1)}$).
\endproclaim

\demo{Proof}
We write $E^k$ as a sequence  $t_\nu^k$ and let $J^k_\mu$ consist
of those $t\in E^k$ where $2^{k-\mu}\le t_{\nu}^k-t_{\nu+1}^k< 2^{k-\mu+1}$
(assuming without loss of generality that the $t^k_\nu$ are decreasing
 in $\nu$).
We clearly have $\card J^k_\mu\lc N(E^k,2^{k-\mu})$.

Let $a^k_\mu$ and $b^k_{\mu}$ denote the endpoints of the equally spaced set
$J^k_\mu$.
Let $\cD^k=\cup_\mu\{a^k_\mu, b^k_\mu\}$, the set of endpoints.
The assertion is implied by the estimate
$$N(\cD^k, 2^{k-j})\lc 2^{j/(1+\beta)} .
\tag 8.2
$$
Let
$L=L_j$ be the smallest integer $\ge 2^{j/(1+\beta)}$.
Note that the set
$\cup_{\mu\ge L_j} J^k_\mu$ is contained in an interval of length
$$
\lc N(E^k, 2^{k-L_j}) 2^{k-L_j}\lc 2^{L_j}(1+L_j)^{-\beta}
2^{k-L_j}
\lc (1+L_j)^{-\beta} 2^k.
$$
This interval can be
 covered by intervals of length $2^{k-j}$ and we need at most
$
(1+L_j)^{-\beta} 2^j$ such intervals to do this. But
$(1+L_j)^{-\beta} 2^j \lc 2^{j/(1+\beta)}$.

We still need to cover the points in $\cD^k$ which do not belong to
$\cup_{\mu\ge L_j} J^k_\mu$. But $\cD^k$ consists just of the $a_\mu^k$ and the
$b_\mu^k$ and there are at most $ 2L_j \lc 2^{j/(1+\beta)}$ points in
$\cD^k$ which are not yet covered. This implies (8.2).

In order to verify the condition (1.11) it suffices  to show
$$\sum_{\mu>n}2^{-\mu} \card(J^k_\mu)\lc 2^{-n} N(E^k, 2^{k-n}) .
\tag 8.3
$$
But if $a^k=\inf_\mu a^k_\mu=\inf E^k$ then the left side of
(8.3) is $\approx 2^{-k}(b^k_n-a^k)$. Moreover  every subinterval 
of length $2^{k-n}$ of $[a^k,b^k_n]$ contains points in $E^k$ and 
therefore
$b^k_n-a^k\lc 2^{k-n} N(E^k, 2^{k-n})$; 
thus (8.3) holds.\qed
\enddemo

\demo{\bf Proof of Theorem 1.1} The set $E^k=\{ 2^k(1+\nu^{-\alpha}):\nu\in \bbZ^+\}$ satisfies
$N(E^k, 2^k\delta)\lc \delta^{1/(\alpha+1)}$ and assertion (i) 
follows 
from 
 Lemma 8.1.1 and Theorem II. On the other hand, the set
$E^k = \{ 2^k (1 + \log^{-\beta} (2+\nu) : \nu \in \bbZ^+ \}$ satisfies
$N(E^k, 2^k \delta) \lc \delta^{-1} [\log(1/\delta)]^{-\beta}$
and assertion (ii) follows from Lemma 8.1.1 and Theorem IV. \qed

\enddemo

\

\subheading{8.2. A counter-example to $\boldkey L^{\boldkey p}$ boundedness for a related 
  maximal function}

Let $E_0$ be any set in $[1,2]$ and define the modified maximal function
$$ \widetilde M_{E_0}f(x) := \sup_{r \in E_0} f * d\sigma(x + re_1)$$
in $\Bbb R^d$, where $e_1$ is a unit vector. If $E_0$
satisfies the regularity assumption $(\cR_p)$, $p< d/(d-1)$ then   the condition 
$\cC_{p,p}$ 
is necessary and sufficient for $L^p$ boundedness 
of $\widetilde M_{E_0}$; indeed a notational modification of the 
proof of Theorem I applies to show the sufficiency.
Note that for sets $E_0$ supported in $[1,2]$ the conditions $\cC_{p,p}$ and $\cC_{p,\infty}$ both amount to the inequality 
$N(E_0,\delta)\lc \delta^{-(d-1)(p-1)}$.
However $L^p$ boundedness and indeed the weak type $(p,p)$ property may
 fail if we drop the regularity assumption.

Let $E_0$ be the 
 middle-halves Cantor set  consisting of all $t=1+ \sum_{j=1}^\infty b_j 4^{-j}$
where $b_j\in \{0,2\}$. Then the Minkowski dimension of $E_0$ is $1/2$ and 
$\widetilde M_{E_0}$ is bounded on $L^p(\Bbb R^2)$ for $p>3/2$ and unbounded on $L^p(\Bbb R^2)$ for $p<3/2$. Moreover  $\cC_{p,p}$ holds for $p=3/2$. We show that nevertheless 
$\widetilde M_{E_0}$  is not of 
weak type $(3/2,3/2)$.

Let $N$ be  large and define 
$$ f(x) := \sum_{i = 1}^N 4^i \chi_{2Ce_1 + B(0,a 4^{-i})}(x),$$
where $C$ is the 
Cantor set $C= \{ \sum_j c_j 4^{-j}: j = 0,1 \}$ and $a$ is small.
Note that $\|f\|_{3/2}\lc N^{2/3}$ (each $i$ contributes an $L^{3/2}$ 
norm of $O(1)$, and the contributions are mostly disjointly supported).

Now  $E_0 + C$ fills out the interval [1,2] and thus the maximal function $\widetilde M_{E_0} f$ 
has size about $N$ on a fixed portion of the unit annulus, thus
$\|\widetilde M_{E_0} f\|_{L^{3/2,\infty}}\ge c N$. This shows that $\widetilde M_{E_0}$ is not 
of weak type $(3/2,3/2)$. 
A closer examination shows that $f$ belongs to the Lorentz space $L^{3/2,s}$ with norm 
$O(N^s)$ so that 
$\widetilde M_{E_0}$ fails to map the  Lorentz space $L^{3/2,s}$ to $L^{3/2,\infty}$ when $s>1$.
Unfortunately this example is too rigid in order to apply to the maximal operator $\cM_{E_0}$
considered in this paper.

\subheading{8.3. Failure of  restricted weak type (2,2) in two dimensions}
We shall now turn to the counterexample mentioned in the introduction and give a proof of
Proposition 1.5.

Suppose that there is a large constant $B$ so that there exists $k$ and $n\ge 100$
such that 
$$N(E^k, 2^{k-2n})\ge B 2^{2n}/n.$$
We then show that
$\|\cM_E\|_{L^{2,1}\to L^{2,\infty}}\ge c\sqrt B$ for some absolute constant $c$.
By rescaling we may assume $k=0$ and $n\gg 1$.

We use the construction of a  Kakeya set as given by Keich \cite{6}, rescaled to a square of 
sidelength $2^{-n}$. It gives us  $\approx 2^n$ rectangles $R_l$  
with sidelengths $2^{-n-3}$ and 
$2^{-2n-6}$ so that $R_l\subset [-2^{-n},2^{-n}]^2$ and the longer side of $R_l$ is parallel
to $e_l:=(\cos l 2^{-n}, \sin l 2^{-n})$, and the union $A=\cup R_l$ has measure 
$\lc 2^{-2n} n^{-1}$. Thus $\|\chi_A\|_{L^{2,1}}
\approx \|\chi_A\|_{2}
\lc 2^{-n} n^{-1/2}$.

Let  $\{I_\nu\}_{\nu=1}^N$ be a cover of the set $E^0$ by dyadic 
intervals of length $2^{-2n}$,
with disjoint interior
so that $N\ge B 2^{2n}/n$.
Let $I_\nu=[a_\nu,b_\nu]$,
and assume $a_\nu<a_{\nu+1}$. We then pick every 
tenth interval $=I_{10 \nu}$, moreover we pick every tenth rectangle
$R_{10 l}$ in the above Kakeya construction.
Let 
$e_l^\perp:=(-\sin l 2^{-n}, \cos l 2^{-n})$ and let 
$R_{l,\nu}$ be the translate $a_{10 \nu} e_{10 l}^\perp+R_{10 l}$. Then the rectangles
$R_{l,\nu}$ are disjoint, 
however on a tenth fraction of each of these rectangles
we have that $\cM_E \chi_A(x)> c2^{-n}$. There are $\approx N 2^n/100$ such rectangles
and thus
$$
\meas\big(\{x: \cM_E\chi_A(x)> c 2^{-n}\}\big)
\ge c'N 2^n 2^{-3n}\gc Bn^{-1} ;
$$
but on the other hand
$\|\chi_A\|_{2}^2/(2^{-2n})\lc n^{-1}$ so that the $L^{2,1}\to L^{2,\infty}$  operator norm 
is  $\gc \sqrt B$. This proves the proposition.\qed

\end
\Refs

\ref \no 1\by J. Bourgain
\paper Estimations de certaines
fonctions maximales
\jour C. R. Acad. Sc. Paris, \vol 310\yr 1985
\endref

\ref  \no 2 \bysame
\paper Averages in the plane  over convex curves
and maximal operators
\jour Jour. Anal. \vol 47 \yr 1986\pages 69--85
\endref



\ref \no  3\by S.Y.A. Chang and R. Fefferman
\paper The Calder\'on-Zygmund decomposition on product domains
\jour Amer. J. Math \vol 104 \yr 1982 \pages 445--468 \endref

\ref \no  4\by M. Christ    \paper Weak type (1,1) bounds for rough
operators \jour Annals of Math. \vol 128 \yr 1988 \pages 19--42
\endref

\ref\no 5\by G. Hardy, J.E. Littlewood and G. P\'olya
\book Inequalities 
\publ Cambridge Univ. Press \publaddr Cambridge\yr 1952
\endref


\ref\no 6\by U. Keich\paper On $L^p$ bounds for Kakeya maximal functions and the Minkowski dimension in $\Bbb R^2$\jour Bull. London Math. Soc\vol 31\yr 1999\pages 213--221\endref


\ref\no 7\by M. Leckband\paper A note on the spherical maximal operator
for radial functions\jour Proc. Amer. Math. Soc.\vol 100\yr 1987
\pages 635--640\endref


\ref\no 8 \by J. Peetre \paper On spaces of Triebel-Lizorkin type
\jour Ark. Mat. \vol 13 \yr 1975 \pages 123--130 \endref


\ref\no 9\by A. Seeger
\paper Remarks on singular convolution operators
\jour  Studia Math. \vol 97 \yr 1990\pages 91--114
\endref

\ref\no 10\by A. Seeger and T. Tao
\paper Sharp Lorentz space estimates for rough operators\jour Math. Annalen \vol 320
\yr 2001\pages
381--415
\endref

\ref\no 11\by A. Seeger, T. Tao and J. Wright
\paper Singular maximal functions and Radon transforms near $L^1$
 \jour preprint
\endref


\ref\no 12\by A. Seeger, S. Wainger and J. Wright\paper
 Pointwise convergence of spherical means
\jour Math. Proc.Cambr. Phil. Soc.\vol 118\yr 1995\endref

\ref\no 13\bysame
\paper Spherical maximal operators on radial functions
 \jour Math. Nachr.\vol 187\yr 1997\pages 95--105\endref


\ref\no 14\by E. M. Stein  \paper Maximal functions:
spherical means \jour Proc. Nat. Acad. Sci.
\vol 73 \yr 1976
\pages 2174--2175
\endref

\ref\no 15\bysame\book Harmonic analysis: Real variable methods,
orthogonality and
 oscillatory integrals\publ Princeton Univ. Press \publaddr Princeton\yr 1993
\endref

\ref\no 16\by T. Wolff\paper Recent work connected with the Kakeya problem
\inbook Prospects in Mathematics (Princeton, N.J., 1996)\pages 129--162\publ Amer. Math. Soc.
\publaddr Providence, R.I.\yr 1999\endref

\endRefs
